\documentclass[11pt]{article}
\usepackage[margin=1in]{geometry} 
\usepackage{amssymb,amsfonts,amsmath,bbm,mathrsfs,stmaryrd,mathtools}
\usepackage{xcolor}
\usepackage{url}

\usepackage{enumerate}

\usepackage{hyperref}
\hypersetup{colorlinks,
             linkcolor=black!75!red,
             citecolor=blue,
             pdftitle={},
             pdfproducer={pdfLaTeX},
             pdfpagemode=None,
             bookmarksopen=true
             bookmarksnumbered=true}

\usepackage{tikz}
\usetikzlibrary{arrows,calc,decorations.pathreplacing,decorations.markings,intersections,shapes.geometric,through,fit,shapes.symbols,positioning,decorations.pathmorphing}

\usepackage{braket}

\usepackage[amsmath,thmmarks,hyperref]{ntheorem}
\usepackage{cleveref}

\creflabelformat{enumi}{#2(#1)#3}

\crefname{section}{Section}{Sections}
\crefformat{section}{#2Section~#1#3} 
\Crefformat{section}{#2Section~#1#3} 

\crefname{subsection}{\S}{\S\S}
\AtBeginDocument{%
  \crefformat{subsection}{#2\S#1#3}%
  \Crefformat{subsection}{#2\S#1#3}%
}

\crefname{subsubsection}{\S}{\S\S}
\AtBeginDocument{%
  \crefformat{subsubsection}{#2\S#1#3}%
  \Crefformat{subsubsection}{#2\S#1#3}%
}

\theoremstyle{plain}

\newtheorem{lemma}{Lemma}[section]
\newtheorem{proposition}[lemma]{Proposition}
\newtheorem{corollary}[lemma]{Corollary}
\newtheorem{theorem}[lemma]{Theorem}

\theoremstyle{nonumberplain}
\newtheorem{theoremN}{Theorem}

\theoremstyle{plain}
\theorembodyfont{\upshape}
\theoremsymbol{\ensuremath{\blacklozenge}}

\newtheorem{definition}[lemma]{Definition}
\newtheorem{example}[lemma]{Example}

\newtheorem{remark}[lemma]{Remark}
\newtheorem{remarks}[lemma]{Remarks}

\crefname{definition}{definition}{definitions}
\crefformat{definition}{#2definition~#1#3} 
\Crefformat{definition}{#2Definition~#1#3} 

\crefname{ex}{example}{examples}
\crefformat{example}{#2example~#1#3} 
\Crefformat{example}{#2Example~#1#3} 

\crefname{exs}{examples}{examples}
\crefformat{examples}{#2example~#1#3} 
\Crefformat{examples}{#2Example~#1#3} 

\crefname{remark}{remark}{remarks}
\crefformat{remark}{#2remark~#1#3} 
\Crefformat{remark}{#2Remark~#1#3} 

\crefname{remarks}{remarks}{remarks}
\crefformat{remarks}{#2remark~#1#3} 
\Crefformat{remarks}{#2Remark~#1#3} 

\crefname{convention}{convention}{conventions}
\crefformat{convention}{#2convention~#1#3} 
\Crefformat{convention}{#2Convention~#1#3} 

\crefname{notation}{notation}{notations}
\crefformat{notation}{#2notation~#1#3} 
\Crefformat{notation}{#2Notation~#1#3} 

\crefname{table}{table}{tables}
\crefformat{table}{#2table~#1#3} 
\Crefformat{table}{#2Table~#1#3}

\crefname{lemma}{lemma}{lemmas}
\crefformat{lemma}{#2lemma~#1#3} 
\Crefformat{lemma}{#2Lemma~#1#3} 

\crefname{proposition}{proposition}{propositions}
\crefformat{proposition}{#2proposition~#1#3} 
\Crefformat{proposition}{#2Proposition~#1#3} 

\crefname{corollary}{corollary}{corollaries}
\crefformat{corollary}{#2corollary~#1#3} 
\Crefformat{corollary}{#2Corollary~#1#3} 

\crefname{theorem}{theorem}{theorems}
\crefformat{theorem}{#2theorem~#1#3} 
\Crefformat{theorem}{#2Theorem~#1#3} 

\crefname{enumi}{}{}
\crefformat{enumi}{(#2#1#3)}
\Crefformat{enumi}{(#2#1#3)}

\crefname{assumption}{assumption}{Assumptions}
\crefformat{assumption}{#2assumption~#1#3} 
\Crefformat{assumption}{#2Assumption~#1#3} 

\crefname{equation}{}{}
\crefformat{equation}{(#2#1#3)} 
\Crefformat{equation}{(#2#1#3)}

\numberwithin{equation}{section}

\theoremstyle{nonumberplain}
\theoremsymbol{\ensuremath{\blacksquare}}

\newtheorem{proof}{Proof}
\newcommand\pf[1]{\newtheorem{#1}{Proof of \Cref{#1}}}
\newcommand\pff[3]{\newtheorem{#1}{Proof of \Cref{#2}#3}}

\newcommand\bC{{\mathbb C}}

\newcommand\bG{{\mathbb G}}

\newcommand\bR{{\mathbb R}}
\newcommand\bS{{\mathbb S}}
\newcommand\bT{{\mathbb T}}

\newcommand\bZ{{\mathbb Z}}

\newcommand\cH{{\mathcal H}}

\newcommand\cQ{{\mathcal Q}}

\newcommand\cS{{\mathcal S}}
\newcommand\cT{{\mathcal T}}

\newcommand\fg{{\mathfrak g}}
\newcommand\fh{{\mathfrak h}}

\newcommand\fk{{\mathfrak k}}
\newcommand\fl{{\mathfrak l}}

\newcommand\ft{{\mathfrak t}}
\newcommand\fu{{\mathfrak u}}

\newcommand\fz{{\mathfrak z}}

\DeclareMathOperator{\id}{id}
\DeclareMathOperator{\im}{im}
\DeclareMathOperator{\ad}{ad}
\DeclareMathOperator{\Ad}{Ad}
\DeclareMathOperator{\End}{\mathrm{End}}
\DeclareMathOperator{\Hom}{\mathrm{Hom}}
\DeclareMathOperator{\Aut}{\mathrm{Aut}}

\newcommand{\cat}[1]{\textsc{#1}}

\newcommand{\qedhere}{\mbox{}\hfill\ensuremath{\blacksquare}}

\title{Subgroup proximity in Banach Lie groups}
\author{Alexandru Chirvasitu}

\begin{document}

\date{}

\newcommand{\Addresses}{{
  \bigskip
  \footnotesize

  \textsc{Department of Mathematics, University at Buffalo}
  \par\nopagebreak
  \textsc{Buffalo, NY 14260-2900, USA}  
  \par\nopagebreak
  \textit{E-mail address}: \texttt{achirvas@buffalo.edu}


}}

\maketitle

\begin{abstract}
  Let $U$ be a Banach Lie group and $G\le U$ a compact subgroup. We show that closed Lie subgroups of $U$ contained in sufficiently small neighborhoods $V\supseteq G$ are compact, and conjugate to subgroups of $G$ by elements close to $1\in U$; this generalizes a well-known result of Montgomery and Zippin's from finite- to infinite-dimensional Lie groups. Along the way, we also prove an approximate counterpart to Jordan's theorem on finite subgroups of general linear groups: finite subgroups of $U$ contained in sufficiently small neighborhoods $V\supseteq G$ have normal abelian subgroups of index bounded in terms of $G\le U$ alone.

  Additionally, various spaces of compact subgroups of $U$, equipped with the Hausdorff metric attached to a complete metric on $U$, are shown to be analytic Banach manifolds; this is the case for both (a) compact groups of a given, fixed dimension, or (b) compact (possibly disconnected) semisimple subgroups. Finally, we also prove that the operation of taking the centralizer (or normalizer) of a compact subgroup of $U$ is continuous (respectively upper semicontinuous) in the appropriate sense.
\end{abstract}

\noindent {\em Key words: Banach Lie group; Banach Lie algebra; analytic manifold; exponential; compact group; Lie subgroup; Hausdorff metric; jet; centralizer; normalizer; semisimple; Hausdorff distance}

\vspace{.5cm}

\noindent{MSC 2020: 22E65; 17B65; 58B25; 32K05; 22E40; 22E30; 22E15; 54E35; 54E50; 47B01}

\tableofcontents

\section*{Introduction}

The ambient Lie groups $U$ of the paper are possibly infinite-dimensional, modeled on Banach spaces: {\it Banach Lie groups}, for short (\cite[\S VI.5]{lang-fund}, \cite[Definition IV.1]{neeb-inf}). The common theme is that of understanding, for a {\it compact} (automatically Lie) subgroup $G\le U$, what ``closeness'' to $G$ entails for other subgroups of $U$.

Although there are various flavors of it (\Cref{def:tpsbgp}), for the purpose of the present less-than-ideally-formal sketch the reader can interpret `close' as `contained in a sufficiently small neighborhood of $G$'. The original impetus and main motivator for the sequel is a result due to Montgomery and Zippin \cite[Theorem 1 and Corollary]{mz-conj}, to the effect that
\begin{itemize}
\item given a compact subgroup $G\le U$ of a {\it finite-dimensional} Lie group $U$;
\item and a neighborhood $V$ of $1\in U$;
\item subgroups of $U$ contained in a sufficiently small neighborhood of $G$ are conjugate to subgroups of $G$ by elements in $V$.
\end{itemize}
Or, even less formally: subgroups of $U$ almost contained in $G$ are almost conjugate to subgroups of $G$. Apart from \cite{mz-conj}, other variations in the literature include, for instance:
\begin{itemize}
\item an analogue for compact subgroups of the group of diffeomorphisms of a compact (hence finite-dimensional) smooth manifold (\cite[Theorem 1]{gk-semic-normal}, \cite[\S VI.7, Theorem 7.1]{om-inflie-bk});
\item the upper semicontinuity of the group of Riemannian isometries of a fixed smooth manifold (i.e. as the Riemannian structure varies, the groups of isometries can only ``jump up'', containment-wise, up to conjugation; \cite[Theorem 8.1]{eb});
\item any number of semicontinuity-flavored cousins thereof, e.g. for automorphism groups of pseudoconvex complex domains: \cite[Theorems 0.1 and 0.2]{gk-pcnvx} or \cite[Theorem 5.1]{gk-surv} and other sources cited in the latter survey.
\end{itemize}

The direct infinite-dimensional analogue of the Montgomery-Zippin result appears below; amended, \Cref{th:closeconj} and \Cref{th:autocpct} read

\begin{theoremN}
  Let $G\le U$ be a compact subgroup of a Banach Lie group.
  \begin{enumerate}[(1)]
  \item\label{item:iscpct} Lie subgroups of $U$ contained in sufficiently small neighborhoods of $G$ are automatically compact.
  \item For any neighborhood $V\ni 1$ in $U$, Lie subgroups of $U$ contained in sufficiently small neighborhoods of $G$ are conjugate to subgroups of $G$ by elements of $V$.  \qedhere
  \end{enumerate}
\end{theoremN}
The compactness statement (i.e. claim \Cref{item:iscpct}) is a non-issue in the finite-dimensional case: locally Euclidean spaces are locally compact.

The thread meanders in a number of directions; in analyzing the above setup for {\it finite} subgroups of $U$, one is led naturally to an approximate version of sorts of Jordan's theorem on finite subgroups of $GL_d(\bC)$ (\cite[Theorem 36.13]{cr-rep}): they all have abelian subgroups of some index bounded by $d$ alone; compare this to \Cref{th:approxjordan} below, which might be of some independent interest:

\begin{theoremN}
  Let $G\le U$ be a compact subgroup of a Banach Lie group.

  There is an $N\in \bZ_{>0}$ such that all finite subgroups $F\le U$ contained in a sufficiently small neighborhood of $G$ have abelian normal subgroups of index $\le N$.  \qedhere
\end{theoremN}

On the other hand, certain families of subgroups of $U$, topologized appropriately, exhibit rich structure (\Cref{th:sameg,th:allsm}):

\begin{theoremN}
  Let $U$ be a Banach Lie group, equipped with a complete metric inducing its topology.
  \begin{enumerate}[(1)]
  \item The following collections of subgroups of $U$, equipped with the Hausdorff distance \cite[Definition 7.3.1]{bbi} attached to the metric on $U$, are all analytic Banach manifolds:
    \begin{itemize}
    \item the set of compact subgroups of any fixed dimension;
    \item the set of compact subgroups of any fixed dimension with a fixed number of connected components;
    \item the set of semisimple (possibly disconnected) compact subgroups. 
    \end{itemize}
  \item For any fixed compact Lie group $G$, the space of continuous embeddings $G\le U$ is a principal $\Aut(G)$-fibration over the Hausdorff-distance-metrized base space of subgroups of $U$ isomorphic to $G$.  \qedhere
  \end{enumerate}
\end{theoremN}

This is this a natural follow-up to some of the work in \cite{cp-us}, focused among other things on spaces of {\it morphisms} $G\to U$ with fixed compact domain (rather than subgroups of $U$), and the result makes use of that earlier material (e.g. \cite[Corollary 2.4]{cp-us}, on the analytic-manifold structure on the space of morphisms $G\to U$).

Finally, in a different direction but broadly on the same theme of subgroup closeness, there is the following result on the (semi)continuity of the centralizer and normalizer operations applied to compact subgroups of a Banach Lie group (\Cref{th:cetnorm}). In the statement, {\it ad-bounded} subsets of $U$ (\Cref{def:adbdd}) are those inducing uniformly bounded adjoint operators on the Lie algebra of $U$. 

\begin{theoremN}
  Let $U$ be a Banach Lie group, equipped with a complete metric inducing its topology, $G\le U$ a compact subgroup, and $B\subset U$ an ad-bounded subset.

  If $K\le U$ is a compact subgroup sufficiently Hausdorff-close to $G$, then
  \begin{enumerate}[(1)]
  \item the intersections of
    \begin{equation*}
      u\left(\text{centralizer of $K$}\right)u^{-1}
      \quad\text{and}\quad
      \left(\text{centralizer of $G$}\right)
    \end{equation*}
    with $B$ are equal for $u\in U$ close to 1; 
  \item and similarly,
    \begin{equation*}
      u\left(\text{normalizer of $K$}\right)u^{-1}\cap B
      \subseteq
      \left(\text{normalizer of $G$}\right)\cap B
    \end{equation*}
    for $u\in U$ close to 1.  \qedhere
  \end{enumerate}
\end{theoremN}
Plainer (but more loosely): the operation of taking centralizers is Hausdorff-metric-continuous up to conjugacy along ad-bounded subsets, and that of taking normalizers is upper semicontinuous in the same sense.

\subsection*{Acknowledgements}

This work was partially supported by NSF grant DMS-2001128.

\section{Preliminaries}\label{se.prel}

The phrase `Lie group', unqualified, always means `analytic real Banach Lie group' (complex Banach Lie groups might appear occasionally, but the reader will be warned). These are the objects that \cite[Chapter III]{bourb-lie-13} specializes to, for instance, when
\begin{itemize}
\item the ground field is either $\bR$;
\item and the ``regularity class'' $C^r$ is $C^{\omega}$ (meaning `analytic'): see \cite[introductory discussion preceding \S III.1]{bourb-lie-13} and \cite[\S 3.2.1]{bourb-vars-17}, which provides the notation.
\end{itemize}

These are also Lie groups of \cite[Definition IV.1]{neeb-inf}, and any number of other sources we will occasionally cite. Although loc.cit. requires only {\it smoothness}, this is no restriction: smooth Banach Lie groups are automatically analytic \cite[Corollary IV.1.10]{neeb-lc}.

Lie {\it subgroups} $G\le U$ are as in \cite[\S III.1.3, Definition 3]{bourb-lie-13}: not only are they automatically closed \cite[\S III.1.3, Proposition 5]{bourb-lie-13}, but the resulting embeddings $\iota:T_1G\to T_1U$ of tangent spaces are required to {\it split}:
\begin{equation*}
  \exists\text{ Banach-space morphism } \pi:T_1U\to T_1G,\quad \pi\circ\iota = \id
\end{equation*}
This splitting requirement is implicit in the notion of `submanifold' employed in the definition and borrowed from \cite[\S 5.8.1]{bourb-vars-17}. On the few occasions when we want to consider closed subgroups that do not meet this splitting requirement, we follow \cite[\S III.1.3, ending Remark]{bourb-lie-13} in referring to these as {\it Lie quasi-subgroups}.

We will also always take it for granted that the Lie algebra $\fu$ of a Lie group has been equipped with a complete norm making it into a {\it Banach Lie algebra} \cite[Definition IV.1]{neeb-inf}: the Lie bracket is continuous, i.e.
\begin{equation}\label{eq:bliealg}
  \|[x,y]\|\le C\|x\|\|y\|,\ \forall x,y\in \fu
  \text{ for some }C>0.
\end{equation}

A smattering of group-theoretic notation:
\begin{itemize}
\item $G_0$ is the connected component containing $1\in G$ (the {\it identity component} of $G$).
\item $Z(G)$ is the center of $G$, and $Z_0(G) = Z_0(G)$ the identity component of that center (i.e. the `0' subscript adorns `$Z$' rather than `$Z(G)$' as a whole).
\item For a Lie group $U$ and $u\in U$ we write $\Ad_u$ for both the conjugation operator
  \begin{equation*}
    u\cdot u^{-1}\in \Aut(U)
  \end{equation*}
  and the operator on the Lie algebra $\fu:=Lie(U)$ induced by $u$-conjugation (i.e. the {\it adjoint action} of $u$ on $\fu$ \cite[\S III.3.12]{bourb-lie-13}). It will not be difficult to distinguish between the two by context.
  
\item For a subgroup $G\le U$,
  \begin{equation*}
    N_U(G):=\{u\in U\ |\ \Ad_u(G)=G\}
  \end{equation*}
  is the {\it normalizer} of $G$ in $U$.
\item Similarly,
  \begin{equation*}
    Z_U(G):=\{u\in U\ |\ ug=gu,\ \forall g\in G\}
  \end{equation*}
  is the {\it centralizer} of $G$ in $U$.
\end{itemize}

\section{Pairs of mutually-close compact subgroups}\label{se:0close}

\cite[Theorem 1 and Corollary]{mz-conj} say, essentially, that in a finite-dimensional Lie group mutually close compact groups are mutual conjugates by elements close to the identity. We are interested in analogues in the context of {\it infinite}-dimensional (Banach) Lie groups.

There are several sensible ways to topologize the collection of (typically Lie, typically compact) subgroups of a Banach Lie group $U$:

\begin{definition}\label{def:tpsbgp}
  Let $U$ be a topological group and $\cS$ a collection of subgroups therein.
  \begin{enumerate}[(1)]

  \item The {\it weak 0-topology}  $\cT_{0,w}$ on $\cS$ has the sets
    \begin{equation*}
      \{\text{subgroups }G\le U\ |\ G\subset V\},\quad V\subseteq U\text{ open}.
    \end{equation*}
    as a basis.

  \item The {\it (strong, or just plain) 0-topology} $\cT_0$ on $\cS$ is defined so that for any $G\in \cS$, a local basis of open neighborhoods of $G$ is given by
    \begin{equation*}
      \{H\in \cS\ |\ H\le G\cdot V\ \text{and}\ G\le H\cdot V\},\quad V=V^{-1}\text{ an open neighborhood of 1}\in U.
    \end{equation*}

  \item Suppose $\cS$ consists of Lie subgroups in the sense of \cite[\S III.1.3, Definition 3]{bourb-lie-13}. 

    The {\it 1-topology} $\cT_1$ on $\cS$ is the subspace topology induced on $\cS$ by the embedding
    \begin{equation*}
      \cS\ni G\xmapsto{\quad} (G,T_1G)\in (\cS,\cT_0)\times \bG(T_1U),
    \end{equation*}
    where
    \begin{itemize}
    \item $\bG(T_1U)$ is the {\it Grassmannian} of $T_1U$, consisting of the latter's complemented closed subspaces, with its standard Banach-manifold topology (\cite[discussion preceding Lemma 3.12]{upm-ban} or \cite[\S 5.2.6]{bourb-vars-17});
    \item and $(\cS,\cT_0)$ denotes $\cS$ equipped with the (strong) 0-topology.
    \end{itemize}
  \end{enumerate}
\end{definition}

\begin{remarks}\label{res:2top}
  \begin{enumerate}[(1)]
  \item The `0' and `1' prefixes are meant as reminiscent of the language of {\it jets} \cite[\S 12.1]{bourb-vars-8-15}: in general, two entities (functions, curves, etc.) define the same {\it $r$-jet} if they are identical ``to order $r$''. Similarly, the 0-topology measures closeness as actual sets (order 0), whereas the 1-topology measures closeness of the sets underlying the groups {\it and} of their Lie algebras (order 1).
    
  \item The weak 0-topology is of course badly non-Hausdorff, regardless of how well-behaved a class of subgroups of $U$ we equip with it: the trivial subgroup, for instance, is dense in $\cS$.
  \item\label{item:0topcomplete} On the other hand, the 0-topology is much better behaved and is familiar (perhaps in disguise) from any number of other contexts.

    Assume, for instance, that $U$ is a Banach Lie group and $\cS$ consists of {\it all} of its closed subgroups. $\cT_0$ is then completely metrizable \cite[\S III.1.1, Proposition 1]{bourb-lie-13} and, having fixed a complete metric, the strong topology is precisely that induced by the {\it Hausdorff distance} \cite[Definition 7.3.1]{bbi} on subsets of $U$. The collection of {\it closed} subgroups of $U$ is thus completely metrizable \cite[Proposition 7.3.7]{bbi} under the strong topology.

    In general, for a metric space $(X,d)$, we will denote by $(\cH(X),d_{\cH})$ its associated space of bounded closed subsets with the Hausdorff distance.
    
  \item In the preceding discussion, one need not even bother with metrics: the strong topology underlies a {\it uniform structure} \cite[\S III.1.1, D\'efinition 1 and \S III.1.2, D\'efinition 3]{bourb-top-1-4} on $\cS$, and the resulting uniform structure on the set of closed subgroups is {\it complete} in the sense of \cite[\S III.3.3, D\'efinition 3]{bourb-top-1-4}.
    
  \item An analogous observation holds for $\cT_1$: it is complete on the family $\cS$ of finite-dimensional Lie subgroups. To see this, recall from \cite[\S 5.2.6]{bourb-vars-17} that a neighborhood of a complemented subspace $V\le \fu$ in the Lie algebra $\fu:=T_1U$ consists of the graphs $\Gamma_f$ of Banach-space morphisms
    \begin{equation*}
      V\xrightarrow{\quad f\quad} W\text{ for some fixed }W\le U\text{ with }\fu=V\oplus W. 
    \end{equation*}
    The space $V$ is close to $\Gamma_f$ precisely when $f$ is small or, equivalently, when the unit ball of $\Gamma_f$ is close to that of $V$ in the Hausdorff distance on closed subsets of $\fu$.

    The upshot of all of this is that $\cT_1$ is the subspace topology induced by
    \begin{equation*}
      \cS\ni G\xmapsto{\quad} (G,\text{unit ball of }T_1G)\in (\cS,\cT_0)\times (\text{bounded bounded subsets of }\fu),
    \end{equation*}
    and completeness follows that of $\cT_0$ (item \Cref{item:0topcomplete} of the present definition) and that of the second factor \cite[Proposition 7.3.7]{bbi}.

  \item\label{item:1cntbl} Having metrized both the Lie group $U$ and its Lie algebra $\fu$, it makes sense to speak of {\it $\varepsilon$-close} groups (for $\varepsilon>0$) in any of the three topologies (assuming the Lie subgroups finite-dimensional in the case of $\cT_1$): a Lie subgroup $H\le U$ is
    \begin{itemize}
    \item $\cT_{0,w}$-$\varepsilon$-close to $G\le U$ if it is contained in the $\varepsilon$-neighborhood of $G$;
    \item $\cT_0$-$\varepsilon$-close to $G$ if each is contained in the other's $\varepsilon$-neighborhood;
    \item and $\cT_1$-$\varepsilon$-close to $G$ if it is $\cT_0$-$\varepsilon$-close and the unit balls of their Lie algebras are $\varepsilon$-close in $\cH(\fu)$ (notation as in \Cref{item:0topcomplete})                
    \end{itemize}
  \item By the same token, we can probe the various topologies (even the non-Hausdorff, non-metrizable $\cT_{0,w}$) through {\it sequence} convergence (as opposed to more sophisticated devices such as {\it nets} or {\it filters}: \cite[\S\S 11 and 12]{wil-top} or \cite[\S 10]{jos-top}).

    Indeed, per point \Cref{item:1cntbl} above, the topologies are all {\it first-countable} \cite[Example 4.4 (b)]{wil-top}: every point has a countable neighborhood base. As explained in \cite[discussion preceding Proposition 1.2]{frnk1}, they are then all {\it sequential}: the topologies are characterized by sequence convergence alone, in the sense that, for instance, a set $V$ is open if and only if every sequence converging to a point in $V$ is eventually in $V$ (see \cite[Proposition 1.1]{frnk1} for alternative characterizations).
  \end{enumerate}
\end{remarks}

The fundamental distinction between 0- and 1-convergence is easy to illustrate:

\begin{example}\label{ex:torwrap}
  Let $\bT^2:=(\bS^1)^2$ be the 2-dimensional torus, and consider its subgroups
  \begin{equation*}
    G_n:=\{(z,z^n)\ |\ z\in \bS^1\}\subset \bT^2,\quad n\ge 1. 
  \end{equation*}
  The sequence $(G_n)_n$ $\cT_0$-converges to all of $\bT^2$, but the Lie algebras
  \begin{equation*}
    \fg_n:=Lie(G_n)\cong \bR
  \end{equation*}
  are all 1-dimensional, so they will not approach the 2-dimensional Lie algebra $Lie(\bT^2) \cong \bR^2$ in the Grassmannian $\bG(\bR^2)$.

  Furthermore, no subsequence of $(G_n)$ $\cT_1$-converges to a closed {\it subgroup} of $\bT^2$ (so the dimension issue cannot be corrected by shrinking the target group): the tangent spaces $\fg_n$ converge to the Lie algebra of a circle factor
  \begin{equation*}
    \bS^1\cong \{1\}\times \bS^1\subset \bT^2,
  \end{equation*}
  but no subsequence will converge to that circle even in the weak 0-topology. 
\end{example}

In working with the various topologies of \Cref{def:tpsbgp}, we will be concerned almost exclusively with families $\cS$ of {\it compact} subgroups of $U$. Recall that these are automatically finite-dimensional Lie subgroups \cite[Theorem IV.3.16]{neeb-lc}, so we will freely discuss their dimensions and any number of other features reliant on their finite dimensionality.

\subsection{The main theorem and (most of) its proof}\label{subse:mz}

The main result of the present section is the following generalization of \cite[Theorem 1 and Corollary]{mz-conj}, from finite- to possibly infinite-dimensional Lie groups.

\begin{theorem}\label{th:closeconj}
  Let $U$ be a Lie group and $G\le U$ a compact subgroup. There is an open neighborhood $V$ of $G$ and a map
  \begin{equation*}
    \left(\text{compact groups contained in $V$}\right)\ni K
    \xmapsto{\quad}
    u_K\in U,
  \end{equation*}
  $\cT_{0,w}$-continuous at $G$, such that
  \begin{itemize}
  \item $u_G=1$;
  \item and the conjugate $u_K\cdot K\cdot u_K^{-1}$ is contained in $G$ for all $K$.
  \end{itemize}  
\end{theorem}

\begin{remark}\label{re:altform}
  Either one of a couple of alternative formulations of \Cref{th:closeconj} might prove more convenient to handle:
  \begin{enumerate}[(a)]
  \item\label{item:ed} For every $\varepsilon>0$ there is a $\delta>0$ such that compact subgroups $K\le U$ $\cT_{0,w}$-$\delta$-close to $G$ are conjugate to subgroups of $G$ by elements $\varepsilon$-close to $1\in U$.
  \item\label{item:seq} Or again: every 0-weakly convergent sequence
    \begin{equation}\label{eq:gntog}
      G_n\xrightarrow[n]{0\text{ weak}}G
    \end{equation}
    has a subsequence $(G_{n_k})_k$ such that
    \begin{equation*}
      u_k\cdot G_{n_k}\cdot u_k^{-1}\le G,\ \forall k
    \end{equation*}
    for a sequence $u_k\to 1\in U$.
  \end{enumerate}
  The switch in perspective among the formulations might, occasionally, be tacit. 
\end{remark}

The argument entails a bit of a detour. The first observation is that for our purposes, there is not much of a distinction between the weak and plain 0-topologies.

\begin{lemma}\label{le:w00}
  Let $G\le U$ be a compact subgroup of a Lie group. A 0-weakly-convergent sequence \Cref{eq:gntog} has a subsequence 0-converging to some closed subgroup $H\le G$.
\end{lemma}
\begin{proof}
  For each $n$, denote by $\ell_n$ the largest distance between a an element of $G_n$ and one of $G$, and set
  \begin{equation*}
    H_n:=\{g\in G\ |\ d(g,G_n)\le \ell_n\}.
  \end{equation*}
  These are compact subsets of $G$, so some subsequence thereof converges to some compact $H\subseteq G$ in $\cH(G)$ (notation as in \Cref{res:2top} \Cref{item:0topcomplete}), because the latter is compact \cite[theorem 7.3.8]{bbi}. Relabel said subsequence so as to avoid changing the indexing; we then have
  \begin{equation*}
    G_n\xrightarrow[n]{\quad 0\quad}H
  \end{equation*}
  in the (plain) 0-topology by construction, and the fact that $H$ is a closed {\it subgroup} of $G$ follows by continuity. 
\end{proof}

The following device will be used repeatedly (thrice, at least), so it seems sensible to have some handy language for it.

\begin{definition}\label{def:elide}
  Consider a 0-weakly-convergent sequence \Cref{eq:gntog}, and suppose we know \Cref{th:closeconj} holds, in its convergent-sequence form of \Cref{re:altform} \Cref{item:seq}, for a sequence of closed normal subgroups
  \begin{equation*}
    H_n\trianglelefteq G_n,
  \end{equation*}
  all mutually isomorphic and also isomorphic to a closed $H\le G$ to which they 0-converge.

  By the hypothesis (that \Cref{th:closeconj} is applicable to $(H_n)_n$), after passing to a subsequence if necessary, we can find
  \begin{equation*}
    u_n\to 1\in U,\quad u_n\cdot H_n\cdot u_n^{-1}\le H.
  \end{equation*}
  The isomorphism hypothesis then requires equality, and said these conjugations transfer the problem to subgroups
  \begin{equation*}
    u_n\cdot G_n\cdot u_n^{-1}\le N_U(H)
  \end{equation*}
  of the {\it normalizer} of $H$ in $U$. Note that $G$ too can be assumed to normalize $H$, since the $G_n$ may as well converge 0-{\it strongly} to $G$ (\Cref{le:w00}). 

  $N_U(H)\le U$ is a Lie subgroup because $H$ is compact \cite[Proposition 2.26]{cp-us}. We can now replace $U$, $G$ and $G_n$ with $N_U(H)/H$, $G/H$ and $G_n/H$.
  
  We will refer to the process just described here as {\it eliding} the normal subgroups $H_n\trianglelefteq G_n$. It has the effect of reducing the problem at hand from the sequence \Cref{eq:gntog} to
  \begin{equation*}
    G_n/H\xrightarrow{\quad n\quad} G/H. 
  \end{equation*}
  Cognates of the word apply, as needed (performing an {\it elision}, etc.).
\end{definition}

\begin{remark}\label{re:normquasi}
  When $U$ is a {\it linear} Lie group, i.e. appropriately embeddable as a group of invertible elements in a Banach algebra \cite[Definition 5.32]{hm4}, the normalizer of {\it any} Lie subgroup (compact or not) is always a Lie {\it quasi-}subgroup \cite[Proposition 5.54 (i)]{hm4}.
\end{remark}

It will be convenient, in discussing families of compact subgroups $G\le U$ topologized in whatever fashion ($\cT_0$, $\cT_1$, etc.), to have first disposed of tori. Ultimately, the reason why these tend to be more difficult to handle is that a $d$-dimensional torus $\bT^d:=(\bS^1)^d$ has a {\it non-compact} automorphism group
\begin{equation*}
  \Aut(\bT^d)\cong GL_d(\bZ). 
\end{equation*}

\Cref{le:cpctauto} formalizes the observation; it is doubtless well known, but seems somewhat difficult to locate in precisely the desired form. Specifically, the (truly minor) annoyances occasioned by the disconnectedness of the group $G$ tend not be worked out in the sources where automorphism groups are discussed.

We will occasionally mention (possibly disconnected) {\it semisimple} compact Lie groups, meaning those whose Lie algebra is semisimple in the usual sense of having no abelian ideals (\cite[\S I.6.1, Definition 1]{bourb-lie-13} or \cite[Definition 6.3]{hm4}).

\begin{lemma}\label{le:cpctauto}
  For a compact Lie group $G$, the group of automorphisms of $G$ operating trivially on the largest central torus $Z_0(G_0)$ of the identity component $G_0$ is compact in the uniform topology.
\end{lemma}
\begin{proof}
  Denote by $A(G)$ and $S$ respectively the automorphism group in the statement and the $G_0$-central torus $Z_0(G_0)$.

  Any automorphism of $G$ induces one on $G/S$, hence a map
  \begin{equation}\label{eq:agtogs}
    A(G)\le \Aut(G)\to \Aut(G/S).
  \end{equation}
  The identity component $\Aut(G/S)_0$ consists of conjugations by elements in $(G/S)_0$, onto which $G_0$ surjects \cite[Lemma 9.18]{hm4}. Since $G_0$ centralizes $S$, the image of \Cref{eq:agtogs} contains the identity component $A(G/S)_0$ and is thus an open (and closed) subgroup. The remainder of the argument branches.

  \begin{enumerate}[(a)]
  \item\label{item:3} {\bf The kernel of \Cref{eq:agtogs} is compact.} An automorphism of $G$ that induces the identity on both $S$ and $G/S$ must be of the form
    \begin{equation*}
      G\ni g\xmapsto{\quad}\psi(g)g\in G
    \end{equation*}
    for a {\it 1-cocycle} $\psi:G\to S$ \cite[\S 5.1]{ser-gal} for the conjugation action of $G$ on $S$, with $\psi\equiv 1$ on $S$. Because $G_0$ centralizes $S$, $\psi$ descends to a {\it morphism} $G_0/S\to S$, which must factor through a finite quotient because $G_0/S$ is semisimple \cite[Theorem 6.4 (vi)]{hm4}. On the other hand cocycles trivial on $G_0$ descend to cocycles $G/G_0\to S$, which form a compact (abelian) Lie group because $G/G_0$ is finite.
    
    In short, the kernel of \Cref{eq:agtogs} fits into a short exact sequence whose extremal terms are compact, so must be compact.

  \item {\bf $\Aut(G/S)$ is compact.} The argument is very similar in spirit to the preceding one. To keep the notation simple, we will assume $G$ is semisimple (though perhaps disconnected) and work with it in place of $G/S$.

    Because $G_0$ is semisimple and connected, its automorphisms can be regarded as automorphisms of its (semisimple) Lie algebra and hence $\Aut(G_0)$ is compact by \cite[Theorem 6.61]{hm4}. We then have a closed-image morphism
    \begin{equation*}
      \Aut(G)\xrightarrow{\quad} \Aut(G_0)\times \Aut(G/G_0)
    \end{equation*}
    into a compact group, and its kernel can again be recovered as a space of cocycles: every automorphism of the kernel is of the form
    \begin{equation}\label{eq:psicoc}
      G\ni g\xmapsto{\quad\alpha\quad}\psi(g)g\in G,\quad
      \psi|_{G_0}\equiv 1,\quad
      \psi:G\to G_0\text{ a 1-cocycle for conjugation}. 
    \end{equation}
    Because on every coset $gG_0=G_0g$ the automorphism $\alpha$ preserves both the left and right actions of $G_0$, its restriction $\alpha|_{gG_0}$ must be multiplication by a {\it central} element in $G_0$. This means that the cocycles $\psi$ of \Cref{eq:psicoc} take values in $Z(G_0)$, and hence constitute a compact Lie group consisting of cocycles
    \begin{equation*}
      G/G_0\to Z(G_0). 
    \end{equation*}
    We can now conclude as before, in \Cref{item:3}.    
  \end{enumerate}
\end{proof}  

Weaker flavors of Lie subgroups of $U$ are sometimes defined by requiring that they be identifiable with closed Lie subalgebras of $\fu:=Lie(U)$ via the exponential map (e.g. \cite[Proposition IV.5 and Definition IV.6]{neeb-inf} or \cite[Definition 5.32]{hm4}). It will be helpful to streamline the language surrounding this sort of setup:

\begin{definition}\label{def:lnrz}
  Let $U$ be a Lie group, $G\le U$ a Lie subgroup, and $\fg\le \fu$ the resulting Lie-algebra inclusion.

  An open neighborhood (typically a ball) $B_{\fu}\subset \fu$ of the origin {\it linearizes} or {\it flattens} $G$ if
  \begin{itemize}
  \item the exponential map $\exp:\fu\to U$ restricts to an analytic isomorphism of $B_{\fu}$ onto its image $B_U:=\exp(B_{\fu})$;
  \item and furthermore,
    \begin{equation}\label{eq:bgbg}
      B_{\fg}:=B_{\fu}\cap \fg
      \xrightarrow{\quad\exp\quad}
      B_G:=B_U\cap G
    \end{equation}
    is also an analytic isomorphism. 
  \end{itemize}
  For variety, we might also say that $B_U\subset U$ (rather than $B_{\fu}\subset \fu$) linearizes/flattens $G$, or that $G$ is {\it $B_{\fu}$- (or $B_U$-)linearized/flat/flattened}.

  Whenever it makes sense, we write `$\log$' for the inverse of partially-defined exponential maps such as \Cref{eq:bgbg}. 
\end{definition}

We will occasionally have to handle centers of compact subgroup $G\le U$ with extra care; in doing this, it will be of some use to know that flattening plays well with centers:

\begin{lemma}\label{le:centflat}
  Let $U$ be a Lie group with Lie algebra $\fu$. An origin neighborhood $B_{\fu}\subset \fu$ that flattens $U$ also flattens the center $Z(G)$ of any Lie subgroup $G\le U$.
\end{lemma}
\begin{proof}
  The center $Z(G)\le U$ is an {\it analytic or integral} (connected) subgroup \cite[\S III.6.2, Definition 1]{bourb-lie-13}, whose Lie algebra $\fz\le \fu$ consists of the vectors fixed by the adjoint representation of $G$ on $\fu$\cite[Proposition 5.54 (ii)]{hm4}.
  
  As in \Cref{def:lnrz}, set $B_U:=\exp(B_{\fu})$. We have to argue that
  \begin{equation*}
    Z(G)\cap B_U\subseteq \exp(\fz\cap B_{\fu}).
  \end{equation*}
  An element $x\in Z(G)\cap B_U$ is fixed by the adjoint representation of $G$; because the exponential map is equivariant for the adjoint representation and the conjugation action, the line in $\fu$ through $\log(x)$ is fixed by $G$ (pointwise), so is contained in $\fz\cap B_{\fu}$. 
\end{proof}

A few other auxiliary observations will highlight the usefulness of flatness/linearization.

\begin{lemma}\label{le:fltcl}
  Let $U$ be Lie group Lie group with Lie algebra $\fu$, $G\le U$ a compact connected subgroup, and $B_{\fu}\subset \fu$ an open ball that flattens $G$ in the sense of \Cref{def:lnrz}.
  
  Consider also a sequence $G_n\le U$ of $B_{\fu}$-flat compact subgroups converging 0-weakly to $G$. For every cluster point
  \begin{equation*}
    \fh \le \fg:=Lie(G)\text{ of }(\fg_n)_n,\ \fg_n:=Lie(G_n),
  \end{equation*}
  the analytic subgroup $H\le G$ with Lie algebra $\fh$ is closed.
\end{lemma}
\begin{proof}
  Because the finite-dimensional Grassmannian $\bG(G)$ is compact, we can assume (after perhaps passing to a subsequence) that the Lie algebras $\fg_n$ converge to $\fh$. 

  Assuming $H$ is {\it not} closed, for some $v\in \fh$ the one-parameter group
  \begin{equation*}
    \{\exp(tv),\ t\in \bR\}\le H
  \end{equation*}
  is not closed in $G$ \cite[Chapter II, Exercise D.4 (v)]{helg}. It follows that for some $a>0$ the intersection
  \begin{equation*}
    \exp([0,a]v)\cap B_U 
  \end{equation*}
  contains a point off $\exp(B_{\fu}\cap \fh)$. There are vectors $v_n\in \fg_n$ converging to $v$, whence
  \begin{equation*}
    \exp([0,a]v_n)
    \xrightarrow[\quad n\quad]{}
    \exp([0,a]v)
  \end{equation*}
  in $\cH(U)$ (\Cref{res:2top} \Cref{item:0topcomplete}). But then for large $n$,
  \begin{equation*}
    \exp([0,a]v_n)\cap B_U
  \end{equation*}
  will contain elements at least $r$ away from $\exp(B_{\fu}\cap \fg_n)$ for some fixed $r>0$, contradicting the flattening hypothesis.
\end{proof}

\begin{lemma}\label{le:toriconv}
  Let $U$ be a Lie group with Lie algebra $\fu$ and
  \begin{equation*}
    G_n\xrightarrow[n]{0\text{ weak}}G,\quad G\le U\text{ compact}
  \end{equation*}
  a 0-weakly convergent sequence of tori, with $G$ and $G_n$ all $B_{\fu}$-flat for an origin-centered open ball $B_{\fu}\subset \fu$.

  There are
  \begin{itemize}
  \item a subsequence $(G_{n_k})_k$;
  \item and isomorphisms
    \begin{equation*}
      \bT^d\xrightarrow[\cong]{\quad\varphi_k\quad}G_{n_k}
      \text{ for some }d\in \bZ_{>0}
    \end{equation*}
  \end{itemize}
  converging uniformly to a morphism $\varphi:\bT^d\to G$.
\end{lemma}
\begin{proof}
  The 0-weak convergence $G_n\to G$ imposes an upper bound on the dimensions of the Lie algebras $\fg_n:=Lie(G_n)$, and implies that subsequences thereof converge to Lie subalgebras of $\fg:=Lie(G)$. We can thus assume, after passing to a subsequence, that
  \begin{equation*}
    \fg_n\xrightarrow[\quad n\quad]{}\fh\le \fg,\quad \dim~\fh=d=\dim~\fg_n,\ \forall n. 
  \end{equation*}
  The Lie algebra $\fh$ will be abelian by continuity, and its associated analytic subgroup $H\le G$ will be a ($d$-dimensional) torus by \Cref{le:fltcl}. We can also assume, after applying \Cref{le:w00} and passing to a subsequence, that $G_n\to H$ in the strong 0-topology.
  
  For large $n$, there are automorphisms $T_n$ close to the identity implementing isomorphisms
  \begin{equation*}
    \fh\xrightarrow[\quad \cong\quad]{T_n} \fg_n;
  \end{equation*}
  this follows, for instance, from the usual description of the analytic structure on the Grassmannian $\bG(\fu)$ \cite[\S 5.2.6]{bourb-vars-17}, whereby an open neighborhood of $\fh$ consists of the graphs of the operators $\fh\to \fl$ for some fixed decomposition
  \begin{equation*}
    \fu\cong \fh\oplus \fl.
  \end{equation*}    
  The kernel
  \begin{equation*}
    \Lambda_n:=\ker\left(\fg_n\xrightarrow{\quad\exp\quad}G_n\right)\subset \fg_n
  \end{equation*}
  is a lattice, and its pullback $T_n^{-1}\Lambda_n\subset \fh$ will similarly be a lattice. For large $n$, these lattices all enjoy a kind of ``uniform rigidity'': for some small origin neighborhood $V\subset \fh$
  \begin{itemize}
  \item no non-zero element of $T_n^{-1}\Lambda_n$ belongs to $V$, e.g. because $V$ can be chosen small enough to as to ensure that $\exp(V)\subset U$ contains no non-trivial subgroups (\cite[Theorem III.2.3]{neeb-lc}, \cite[Theorem 1.4.2]{om-transf}).
  \item and there is an upper bound on the volume of a fundamental domain of $T_n^{-1}\Lambda_n$, because all such fundamental domains can be covered by $N$ translates of $V$ for some $N\in \bZ_{>0}$ independent of $n$ ($H$ is covered finitely many translates of $\exp(V)$, $G_n\to H$ 0-strongly, etc.).
  \end{itemize}
  {\it Mahler selection} \cite[Theorem 1.1]{mvg-mahl} then shows that some subsequence of $(T_n^{-1}\Lambda_n)_n$ converges to a lattice $\Lambda\subset \fh$, whence by \cite[Theorem 5]{macb-homeo} an embedding
  \begin{equation*}
    \bZ^d\cong \Lambda\subset \fh
  \end{equation*}
  compact-uniformly close to embeddings
  \begin{equation*}
    \bZ^d\cong T_n^{-1}\Lambda_n\subset \fh.
  \end{equation*}
  This gives morphisms
  \begin{equation*}
    \bT^d\cong \fh/\Lambda\to G_n \cong \fh/T_n^{-1}\Lambda_n
  \end{equation*}
  with differentials close to the identity for large $n$, which is what was desired. 
\end{proof}

\begin{remark}\label{re:flatnec}
  Flatness is crucial in \Cref{le:toriconv}: the sequence in \Cref{ex:torwrap} fails to satisfy the conclusion of the lemma precisely because it does not meet the linearization requirement.
\end{remark}

Now that we know (by \Cref{le:toriconv}) that convergent sequences of $B$-flat tori behave well, we complement this by correcting for {\it non-}flatness.

\begin{lemma}\label{le:makeflat}
  Let $U$ be a Lie group with Lie algebra $\fu$ and \Cref{eq:gntog} a 0-weakly-convergent sequence of compact subgroups of $U$. There are 
  \begin{itemize}
  \item an arbitrarily small origin neighborhood $B_{\fu}\subset \fu$;
  \item and a subsequence $(G_{n_k})_k$ whose members can be enlarged to $G_{n_k}\le \widetilde{G_{n_k}}$ so that
    \begin{equation*}
      \widetilde{G_{n_k}}\xrightarrow[k]{0\text{ strong}}H
    \end{equation*}
    for some closed $H\le G$;
  \item and the largest central tori $Z\left(G_{n_k,0}\right)_0$ of the connected components $G_{n_k,0}$ are all $B_{\fu}$-flat. 
  \end{itemize}
\end{lemma}
\begin{proof}
  With every passage to a subsequence, we will simply label that subsequence back to the original indexing to keep the notation light. To the same purpose, set $Z_n:=Z_0(G_{n,0})$.
  
  First, by \Cref{le:w00}, we can always assume that $(G_n)_n$ and $(Z_n)_n$ are both strongly 0-convergent to respective closed subgroups $T\le H$ of $G$. Because $Z_n$ are tori, so is $T$. Furthermore, because the Lie algebras
  \begin{equation*}
    \fz_n:=Lie(Z_n)
  \end{equation*}
  are close to subspaces of $\ft:=Lie(T)$ in the Grassmannian $\bG(\fu)$, the dimensions $\dim Z_n$ are bounded above by $\dim T$ (for large $n$). If $\dim Z_n=\dim T$ for infinitely many $n$, then {\it that} subsequence will be $B$-flattened for a sufficiently small origin-centered ball $B\subset \fu$. There is thus no harm in assuming that $\dim Z_n<\dim T$.

  Consider a ball $B_{\fu}\subset \fu$ that linearizes $G$ and $T$, and suppose it {\it doesn't} linearize any of the $Z_n$ (otherwise we could pass to a subsequence that {\it is} flattened, etc.). Because
  \begin{equation*}
    B_U\cap Z_n\xrightarrow[n]{\quad\text{in $\cH(U)$}\quad} B_U\cap T,\quad B_U:=\exp(B_{\fu})
  \end{equation*}
  and the limit is a strictly higher-dimensional manifold, it is not difficult to see that for large $n$ the set
  \begin{equation*}
    \exp^{-1}\left(B_U\cap Z_n\right)\subset \fz_n
  \end{equation*}
  contains some open subset disjoint from $B_{\fu}\cap \fz_n$. We can then find some finite-order
  \begin{equation*}
    x\in \left(B_U\cap Z_n\right)\setminus \exp\left(B_{\fu}\cap \fz_n\right);
  \end{equation*}
  the one-parameter group
  \begin{equation*}
    {}_x\bS^1:=\{\exp~tx\ |\ t\in \bR\}
  \end{equation*}
  will then be a circle (as the notation suggests), centralized by all of $G_{n,0}$ (because $x$ is so centralized). The finite component group $\pi_0(G_n):=G_n/G_{n,0}$ operates on the centralizer of $G_{n,0}$ by conjugation, thus producing finitely many conjugates of ${}_x\bS^1$.

  Now, ${}_x\bS^1$ itself is a finite union of $G_n$-translates of the segment $\exp([0,\log~x])$; the latter is contained in $B_U$ and $G_n$ is 0-weakly-close to $G$, so ${}_x\bS^1$ will also be 0-weakly-close to $G$. This same argument goes through for all conjugates
  \begin{equation*}
    \Ad_s\left({}_x\bS^1\right),\ s\in \pi_0(G_n)
  \end{equation*}
  so that finally, the groups
  \begin{equation*}
    G_n\cdot \prod_{s\in \pi_0(G_n)}\Ad_s \left({}_x\bS^1\right),\quad n\in \bZ_{>0}
  \end{equation*}
  (with the circle factors all mutually commuting and centralizing $G_{n,0}$) are all 0-weakly-close to $G$.

  This has the effect of enlarging dimensions of the centers of $G_n$ strictly, so the procedure cannot continue indefinitely: eventually, after having been properly enlarged $G_n$ to $\widetilde{G_n}$, the latter must (perhaps after passing to a subsequence) be flat for some small $B\subset \fu$.

  As for a subsequence converging to $H\le G$ strongly, this can always be arranged by an application of \Cref{le:w00}.
\end{proof}

\pff{th:closeconj-red-to-ss}{th:closeconj}{ (reduction to semisimple groups)}
\begin{th:closeconj-red-to-ss}
  Recall that `semisimple' means the Lie algebra is semisimple; it does not entail connectedness. To be more precise on what the current branch of the proof is meant to do, we address the problem in its convergent-sequence form, per \Cref{re:altform} \Cref{item:seq}, and reduce the problem to the case where the individual $G_n$ are semisimple.
  
  If \Cref{eq:gntog} is a 0-weakly-convergent sequence of compact groups, we can assume by \Cref{le:makeflat,le:w00} that the connected centers
  \begin{equation*}
    Z_{n}:=Z_0(G_{n,0})
  \end{equation*}
  are all $B$-flattened by some small origin neighborhood $B\subset \fu$, and that the $G_n$ and $Z_n$ 0-strongly-converge to closed subgroup $Z\le H$ of $G$.
  
  Further, \Cref{le:toriconv} affords us the assumption that $Z_n$ and $Z$ are images of embeddings
  \begin{equation*}
    \bT^d\xrightarrow[\cong]{\quad\varphi_n\quad}Z_n \le G
    \quad\text{and}\quad
    \bT^d\xrightarrow[\cong]{\quad\varphi\quad}Z \le G,
  \end{equation*}
  with $\varphi_n\to \varphi$ in the uniform topology. Now, uniformly close {\it morphisms} $\bT^d\to U$ are conjugate by elements close to $1\in U$ \cite[Corollary 2.6]{cp-us}, so the $Z_n$ can be elided per \Cref{def:elide}. This shifts the problem from \Cref{eq:gntog} to
  \begin{equation*}
    G_n/Z\xrightarrow{\quad n\quad} G/Z,
  \end{equation*}
  with the $G_n/Z$ semisimple by construction. 
\end{th:closeconj-red-to-ss}

This is also a good time to dispose of another particular case:

\pff{th:closeconj-abix}{th:closeconj}{ (abelian subgroups of bounded index)}
\begin{th:closeconj-abix}
  The assumption here is that the groups $G_n$ of \Cref{eq:gntog}, converging 0-weakly to $G$, all have abelian subgroups $H_n\le G_n$ with indices
  \begin{equation*}
    [G_n:H_n] \le N,\ \forall n
  \end{equation*}
  bounded uniformly in $n$.

  There is no harm in assuming the $H_n$ normal (since we can replace each $H_n$ with the intersection of all of its conjugates under $\le N$ representatives for $G_n/H_n$) and, by \Cref{le:w00}, that
  \begin{equation}\label{eq:hgconv}
    G_n\xrightarrow[\quad n\quad]{} G
    \quad\text{and}\quad
    H_n\xrightarrow[\quad n\quad]{} H\le G
  \end{equation}
  0-strongly. $H$ will again be abelian by continuity, and the index $[G:H]$ is also bounded by $N$.

  So long as $\dim H_n<\dim H$ for all sufficiently large $n$, we can enlarge the dimensions of $H_n$ by the same technique as in the proof of \Cref{le:makeflat}. The upshot is that we can assume the identity components $H_{n,0}$ to be tori and
  \begin{equation*}
    \dim H_n = \dim H. 
  \end{equation*}
  The strong convergence \Cref{eq:hgconv} implies that the numbers $|H_n/H_{n,0}|$ of components eventually stabilize to $|H/H_0|$, so that the indices $[G_n:H_{n,0}]$ are also uniformly bounded.
  
  We can now elide the $H_n$, as in the previous argument, making the substitutions
  \begin{equation*}
    U\xmapsto{\quad}N_U(H)/H,
    \quad
    G\xmapsto{\quad} G/H,
    \quad
    G_n\xmapsto{\quad} G_n/H_{n,0}. 
  \end{equation*}
  We have just observed that the latter groups are finite of bounded order, so after possibly passing to some subsequence we can assume that $G_n$ are images of embeddings
  \begin{equation*}
    K\xrightarrow[\quad\cong\quad]{\varphi_n} G_n\le U
  \end{equation*}
  of a single finite group $K$. 

  Because $G_n\to G$ 0-strongly and $G$ is compact, the collection
  \begin{equation*}
    \left\{\varphi_n(x)\ |\ n\in \bZ_{>0}\right\}
  \end{equation*}
  of images of $x$ is relatively compact for every point $x\in K$. Some subsequence of $(\varphi_n)_n$ will then converge uniformly to a morphism $\varphi:K\to U$ by {\it Ascoli's} \cite[Theorem 47.1]{munk}. The image $\varphi(K)$ must be $G$ by continuity, and the conclusion follows from the fact that uniformly close morphisms $\varphi:K\to U$ are conjugate by elements $u\sim 1$ of $U$ \cite[Corollary 2.6]{cp-us}.
\end{th:closeconj-abix}

By contrast to tori, {\it semisimple} compact groups (to which \Cref{th:closeconj} has just been reduced) are much better behaved as far as convergence in the various topologies goes:

\begin{proposition}\label{pr:ssgood}
  Let $G\le U$ be a compact subgroup of a Banach Lie group and \Cref{eq:gntog} a sequence of compact connected semisimple Lie groups 0-weakly converging to $G$.

  There is a compact, connected, semisimple Lie group $K$ with morphisms
  \begin{equation*}
    \left(K\xrightarrow[\cong]{\quad\varphi_k\quad}G_{n_k}\right)
    \xrightarrow[k]{\quad\text{uniform}\quad}
    \left(K\xrightarrow[\cong]{\quad\varphi\quad}H\le G\right)
  \end{equation*}
  for some subsequence $(G_{n_k})_k$.

  Furthermore, $(G_{n_k})_k$ converges to $H$ in the $\cT_1$ topology.
\end{proposition}
\begin{proof}
  As before, we consistently relabel subsequences so as to avoid double indexing. Write $\fu$, $\fg$ and $\fg_n$ for the Lie algebras of $U$, $G$ and $G_n$ respectively.
  
  First, we can assume, by the compactness of the Grassmannian $\bG(\fg)$ that $\fg_n\to \fh$ for some $\fh\le \fg$, which will automatically be a Lie subalgebra by continuity. There is also no loss in assuming all $\fg_n$ are mutually isomorphic:

  The $\fg_n$ all have the same dimension $\dim \fh$ for sufficiently large $n$, and there are finitely many isomorphism classes of complex semisimple Lie algebras in any given dimension (by the well-known classification \cite[\S X.6, Table IV]{helg}), hence also finitely many isomorphism classes of compact semisimple Lie algebras, by the correspondence between the two \cite[Corollary III.7.3]{helg}. It follows that some subsequence of $(\fg_n)_n$, which we can harmlessly substitute for the original, consists of mutually isomorphic Lie algebras.

  There are finitely many isomorphism classes of compact, connected, semisimple Lie groups with a given Lie algebra, so (again, after passing to a subsequence), there will be isomorphisms
  \begin{equation*}
    K\xrightarrow[\cong]{\quad \varphi_n\quad}G_n. 
  \end{equation*}
  The sets 
  \begin{equation*}
    \left\{\varphi_n(x)\ |\ n\in \bZ_{>0}\right\},\ x\in K
  \end{equation*}  
  are relatively compact as in the preceding proof. On the other hand, it follows from the boundedness result of \Cref{le:lalgequicont} that $(\varphi_n)_n$ is {\it equicontinuous} \cite[Definition following Theorem 45.1]{munk}, so another application of Ascoli's theorem ensures uniform convergence  
  \begin{equation*}
    \varphi_n\xrightarrow[\quad n\quad]{}\left(K\xrightarrow[]{\quad\varphi\quad} U\right).
  \end{equation*}  
  The image $\varphi(K)$ must be $H$ because we were already assuming 0-strong convergence $G_n\to H$, and the injectivity of $\varphi$ follows from the injectivity of $\varphi_n$: $\ker\varphi$ would have to be finite in any case, for dimension reasons, $\varphi_n\left(\ker\varphi\right)$ are finite subgroups contained, for large $n$, in arbitrarily small subgroups of $U$, but the latter ``has no small subgroups'' (by \cite[Theorem III.2.3]{neeb-lc}, \cite[Theorem 1.4.2]{om-transf} again).
  
  The uniform convergence $\varphi_n\to \varphi$ now implies that $\varphi_n$ is conjugate to $\varphi$ by $u_n\sim 1$ for large $n$ \cite[Corollary 2.6]{cp-us}, hence the last claim on $\cT_1$-convergence.
\end{proof}

Recall that a {\it compact} Lie algebra \cite[Definition 6.1]{hm4} is a Lie algebra of a compact Lie group. 

\begin{lemma}\label{le:lalgequicont}
  Let $\fk$ be a compact semisimple Lie algebra and $\fu$ a Banach Lie algebra.

  Given a sequence
  \begin{equation*}
    \psi_n:\fk\xrightarrow{\quad} \fu
  \end{equation*}
  of embeddings whose images $\fg_n:=\psi_n(\fk)$ converge in the Grassmannian $\bG(\fu)$, the norms of $\psi_n$ are uniformly bounded above, as are those of the inverses $\psi_n^{-1}:\fg_n\to \fk$.
\end{lemma}
\begin{proof}  
  The bound in question refers to a norm topologizing $\fk$; it can be arbitrary (and fixed throughout the proof), since all norms will induce the same topology.

  The isomorphisms
  \begin{equation*}
    \fk\xrightarrow[\cong]{\quad\psi_n\quad}\fg_n
  \end{equation*}
  (denoted abusively by the same symbols as the original $\psi_n$) will preserve the respective {\it negative Killing forms}
  \begin{equation*}
    \braket{x|y}_{\fl}:= -\text{trace}~\ad_x \ad_y,\quad ad_x:=[x,-]\in \End(\fl)
  \end{equation*}
  of the Lie algebras in question: $\fl=\fk$ or $\fl=\fg_n$. These forms are positive definite \cite[Theorem 6.4 (ix) and Theorem 6.6]{hm4}, so they induce their own respective norms $\|\cdot\|_{\fl,k}$ (the `k' subscript stands for `Killing'). The claim boils down, then, to showing that there are constants $C_i>0$ with
  \begin{equation*}
    C_0\|\cdot\|_{\fg_n,k} \le \|\cdot\|\le C_1 \|\cdot\|_{\fg_n,k},\ \forall n,
  \end{equation*}
  where the middle norm is that inherited by $\fg_n$ from the original ambient Banach Lie algebra $\fu$. Or, in short: the Killing norms are equivalent to the original norm uniformly in $n$.

  It is the {\it lower} bound that is interesting: for the upper (which does not require the full force of the hypothesis), recall that the ambient (Banach) Lie algebra $\fu$ satisfies an inequality of the form \Cref{eq:bliealg}, whence
  \begin{equation*}
    \|\ad_x^2\| = \|[x,[x,-]]\|\le C^2\|x\|^2.
  \end{equation*}
  The $\braket{-|-}_{\fg_n}$-negative-definite operator $\ad_x^2$, $x\in \fg_n$ thus has eigenvalues absolutely bounded by $C^2\|x\|^2$, so
  \begin{equation*}
    \|x\|_{\fg_n,k} = \sqrt{-\mathrm{trace}~\ad_x^2}\le \sqrt{\dim\fg_n}\cdot C \|x\|,\ \forall x\in \fg_n.
  \end{equation*}
  
  On to $C_0$ then. Suppose there {\it isn't} such a bound. This means that after passing to an appropriate subsequence, we can find $x_n\in \fg_n$ with
  \begin{equation*}
    \|x_n\|=1
    \quad\text{and}\quad
    \|x_n\|_{\fg_n,k}\to\infty.
  \end{equation*}
  Because, once more, $\ad_{x_n}^2\in \End(\fg_n)$ is negative definite with respect to $\braket{-|-}_{\fg_n}$, this absolutely-large-trace condition means that it has arbitrarily small negative eigenvalues for increasing $n$. In conclusion, there are
  \begin{equation}\label{eq:xxy}
    y_n\in \fg_n,\ \|y_n\|\to 0
    \quad\text{with}\quad
    [x_n,[x_n,y_n]] = \ad_{x_n}^2(y_n) = \lambda_n y_n\text{ of norm }1. 
  \end{equation}
  
  Recall that we are assuming convergence $\fg_n\to \fg$ in $\bG(\fu)$. After perhaps passing to subsequences, we can thus assume that $x_n$ converge to an element $x\in \fg$ of norm 1, and similarly for
  \begin{equation*}
    \lambda_n y_n\xrightarrow[\quad n\quad]{} z\in \fg. 
  \end{equation*}
  The $y_n$s themselves will approach 0 of course, so in the limit \Cref{eq:xxy} reads
  \begin{equation*}
    0 = [x,[x,0]] = z\text{ of norm }1:\text{ a contradiction}. 
  \end{equation*}
\end{proof}

A ``smallness'' requirement (such as the Grassmannian-convergence hypothesis) on the family of morphisms $\psi_n$ of \Cref{le:lalgequicont} is necessary; the statement is certainly not true without some such restriction:

\begin{example}\label{ex:casimir}
  Let $K:=SU_2$ and $\fk$ its Lie algebra, and take for $U$ the unitary group of an infinite-dimensional Hilbert space. All irreducible unitary representations of $K$ (classified, as usual, by {\it dominant weights}: \cite[\S 21.2, Corollary]{hmph}, \cite[\S 23.2]{fh}, etc.) can be realized in $U$.

  Ranging over all irreducible $K$-representations, the {\it Casimir element} \cite[\S 22.1]{hmph}
  \begin{equation*}
    c\in \text{enveloping algebra of }\fk\otimes_{\bR}\bC
  \end{equation*}
  acts as a non-negative scalar in each, unbounded over the entire family of representations \cite[Exercise 23.4]{hmph}. It follows, then, that the morphisms $\fk\to \fu$ resulting from representations $K\to U$ are not uniformly bounded.
\end{example}

\pff{th:closeconj-red-to-fin}{th:closeconj}{ (reduction to finite groups)}
\begin{th:closeconj-red-to-fin}
  Once more, this means boiling down the problem, in its \Cref{re:altform} \Cref{item:seq} incarnation, to the case where the individual $G_n$ (which 0-converge to $G$) are finite.

  By the earlier arm of the proof, we can assume the $G_n$ are semisimple, in the sense that their Lie algebras $\fg_n$ are. By \Cref{pr:ssgood}, we can assume that the identity components $G_{n,0}$ are images of embeddings
  \begin{equation*}
    \varphi_n:K\to U\quad\text{for a compact, connected, semisimple Lie group }K,
  \end{equation*}
  converging uniformly to an embedding $\varphi:K\to G$. Because mutually close morphisms $K\to U$ are conjugate by unitaries close to $1\in U$ (\cite[Corollary 2.6]{cp-us} again),
  \begin{equation*}
    G_{n,0}\trianglelefteq G_n
  \end{equation*}
  are subject to elision as described in \Cref{def:elide}, so we are reduced to handling the sequence $(G_n/G_{n,0})_n$ of finite groups.  
\end{th:closeconj-red-to-fin}

\subsection{Completing the main proof: an approximate Jordan theorem}\label{subse:jordan}

The title refers to the celebrated result stated on \cite[p.91, Th\'eor\`eme I]{jord} (see also \cite[Theorem 36.13]{cr-rep} or \cite[Theorem 8.29]{ragh} for instance). In a refinement \cite[p.281, statement I]{bw-jordan}, it says that a finite subgroup $F\le U$ of a finite-dimensional connected Lie group $U$ has an abelian normal subgroup
\begin{equation*}
  A\trianglelefteq F
  \quad\text{with}\quad
  |F/A|\le n(K) 
\end{equation*}
for some $n(K)$ depending only on a maximal compact subgroup $K$ of $U$.

In light of the reduction of \Cref{th:closeconj} to a statement about 0-weakly-convergent sequences of finite subgroups of $U$, it is perhaps reasonable to expect that such a result would feature here. `Approximate' comes in because we are interested in finite groups that are {\it close to} (bot not, a priori, contained in) compact subgroups of our possibly infinite-dimensional $U$. In short, then, the statement is:

\begin{theorem}\label{th:approxjordan}
  Let $U$ be a Banach Lie group and $G\le U$ a compact subgroup.

  For sufficiently small neighborhoods $V\supseteq G$, all finite subgroups $F\le U$ contained in $V$ have abelian normal subgroups whose indices are bounded uniformly by a constant depending only on the embedding $G\le U$. 
\end{theorem}

Assuming the statement just momentarily, we can complete the proof of \Cref{th:closeconj}:

\pff{th:closeconj-fin}{th:closeconj}{ (conclusion)}
\begin{th:closeconj-fin}
  Per the preceding reductions, suppose $G_n\to G$ is a 0-weakly convergent sequence with $G$ compact and $G_n$ finite. \Cref{th:approxjordan} places an upper bound on the index of a normal abelian subgroup of $G_n$ (with varying $n$). This case, of $G_n$ having abelian subgroups of bounded index, has already been treated (preceding \Cref{pr:ssgood}), so we are done.
\end{th:closeconj-fin}

Returning to the remaining task:

\pf{th:approxjordan}
\begin{th:approxjordan}
  In what must by now be a familiar fashion, the claim can be rephrased in terms of convergence: if \Cref{eq:gntog} is a 0-weakly-convergent sequence with $G$ compact and $G_n$ finite, then
  \begin{equation*}
    \min\{|G_n/A_n|\ |\ A_n\trianglelefteq G_n\text{ normal abelian}\}
  \end{equation*}
  is bounded in $n$. 

  Assume, in fact, that $G_n\to G$ 0-{\it strongly} (as permitted by \Cref{le:w00}). For each $n$, consider an injection
  \begin{equation*}
    \iota_n:G_n\to G,\ d(\iota x,x)<\varepsilon_n,\ \forall x\in G_n
  \end{equation*}
  for $\varepsilon_n\to 0$. The maps $\iota_n$ will of course not, in general, be {\it group} morphisms, but because they move all elements just a little, they are {\it approximate} group morphisms. In particular, the subsets
  \begin{equation*}
    F_n:=\iota_n G_n\subseteq G,
  \end{equation*}
  equipped with their group structures transported over from $G_n$ via $\iota_n$, constitute {\it $\varepsilon'_n$-approximations} of $G$ in the sense of \cite[p.105]{turing} for some $\varepsilon'_n\to 0$.

  The main result of loc.cit., \cite[Theorem 2]{turing}, ensures that $G$ itself has an abelian subgroup of finite index (so in particular, its identity component is a torus). This is not quite what we need to prove, but the proof of that result can be deconstructed to yield the conclusion desired here.

  Recall how that proof proceeds, via the auxiliary result of \cite[Theorem 1]{turing}:
  \begin{itemize}
  \item Assuming the group $G$ has a faithful $d$-dimensional representation (which compact Lie groups always do \cite[Proposition 5.33 (iv)]{hm4});
  \item and is $\eta$-approximated by a finite group $F$ for some appropriately small $\eta$;
  \item it follows that a quotient of $F$ by a subgroup $N\trianglelefteq F$ contained in an $\varepsilon'$-neighborhood of $G$ also has a faithful $d$-dimensional representation;

  \item and $\varepsilon'$ can be chosen arbitrarily small with $\varepsilon\to 0$.
  \end{itemize}
  In our case, $F$ will be one of the $F_n$. Its normal subgroup $N$ referred to above (again, for the group structure transported via $\iota_n$) will be contained in a small neighborhood of $1\in G$; because $\iota_n$ moves every element only a little, the normal subgroup
  \begin{equation*}
    \iota_n^{-1}N\trianglelefteq G_n
  \end{equation*}
  will also be contained in a small identity neighborhood. By the no-small-subgroups \cite[Theorem III.2.3]{neeb-lc},
  \begin{equation*}
    \iota_n^{-1}N\text{ is trivial }\Rightarrow N\le F_n\text{ is trivial}. 
  \end{equation*}
  The upshot, then, is that for sufficiently large $n$ the finite group $G_n\cong F_n$ has a faithful $d$-dimensional representation; the conclusion now follows from the usual Jordan \cite[Theorem 36.13]{cr-rep}.
\end{th:approxjordan}

\section{1-jet closeness and analytic manifolds of subgroups}\label{se:1close}

By \Cref{th:closeconj}, compact subgroups 0-weakly-close to a compact $G\le U$ can be conjugated into $G$ by elements close to $1\in U$. On the other hand, being conjugates by $u\sim 1$ entails closeness in the $\cT_1$ topology. There is, then, a sense in which 0-closeness can be leveraged into 1-closeness. We formalize this in \Cref{pr:sametop}, clarifying some of the relationships between the various topologies: assuming compactness, $\cT_1$ convergence means 0-convergence plus relatively mild numerical constraints.

\begin{proposition}\label{pr:sametop}
  Let $U$ be a Lie group and $\cS$ a family of compact subgroups thereof. The following topologies on $\cS$ are all equal (with $\bZ_{\ge 0}$ topologized discretely):

  \begin{enumerate}[(a)]

  \item\label{item:1} the 1-topology $\cT_1$;

  \item\label{item:0s} the subspace topology induced by
    \begin{equation*}
      \cS\ni G\xmapsto{\quad} (G,\dim G)\in (\cS,\cT_{0})\times \bZ_{\ge 0}.
    \end{equation*}

  \item\label{item:0w} the subspace topology induced by    
    \begin{equation*}
      \cS\ni G\xmapsto{\quad} (G,\dim G,|G/G_0|)\in (\cS,\cT_{0,w})\times \bZ_{\ge 0}\times \bZ_{\ge 0}
    \end{equation*}
    where $G_0\le G$ is the identity component;

  \end{enumerate}
\end{proposition}
\begin{proof}
  We prove that sequence convergence in either of the topologies implies convergence in the others, in he stated order. The abbreviation `(x) $\Longrightarrow$ (y)' stands for
  \begin{equation*}
    \text{convergence in (x) }\Longrightarrow\text{ convergence in (y)}.
  \end{equation*}
  \begin{enumerate}[]
  \item {\bf \Cref{item:1} $\Longrightarrow$ \Cref{item:0s}.} This is clear: the $\cT_1$ topology is stronger than $\cT_0$ by definition, and the dimension is a continuous function on the subspace of the Grassmannian $\bG:=\bG(Lie(U))$, with
    \begin{equation*}
      \bG_n:=\{V\in \bG\ |\ \dim~V=n\}\subset \bG
    \end{equation*}
    both open and closed \cite[\S 5.2.6]{bourb-vars-17}.

  \item {\bf \Cref{item:0s} $\Longrightarrow$ \Cref{item:0w}.} $\cT_0$ is stronger than $\cT_{0,w}$ by definition, so it remains to argue that convergence
    \begin{equation}\label{eq:genericonv}
      G_n\xrightarrow[\quad n\quad]{}G
    \end{equation}
    in \Cref{item:0s} means convergence in component-number (i.e. the stabilization of that number).

    We may as well assume $G_n$ and $G$ all have the same dimension. \Cref{th:closeconj} ensures that for sufficiently large $n$ the group $G_n$ is conjugate to a subgroup of $G$; equidimensional {\it connected} compact Lie groups cannot contain each other properly, so a conjugation
    \begin{equation*}
      G_n\xrightarrow[]{\quad \Ad_u\quad} G,\ u\in U
    \end{equation*}
    will map every component of $G_n$ onto a component of $G$. Since furthermore $G_n\to G$ 0-{\it strongly}, {\it all} components of $G$ must be obtained in this manner for sufficiently large $n$.

    In summary: for large $n$, $G_n$ and $G$ all have the same number of components.

  \item {\bf \Cref{item:0w} $\Longrightarrow$ \Cref{item:1}.} The claim here is that 0-weak convergence \Cref{eq:genericonv} together with stabilization of dimensions and component-numbers entails 1-convergence. The argument is very similar to the preceding one:

    By \Cref{th:closeconj} sufficiently-high-index $G_n$ are conjugate to subgroups of $G$ by elements $U\ni u\sim 1$. Because $G_n$ and $G$ all have the same dimension, conjugation identifies the components of $G_n$ with some of the components of $G$. But we are also assuming (eventual, for large $n$) equality of component numbers, so in fact
    \begin{equation*}
      u_n\cdot G_n\cdot u_n^{-1} = G,\ u_n\sim 1.
    \end{equation*}
    This, naturally, also entails closeness of tangent spaces in $\bG(Lie(U))$. 
  \end{enumerate}
  This concludes the proof. 
\end{proof}

An immediate consequence:

\begin{corollary}\label{cor:samedimsamecomp}
  Let $U$ be a Lie group and $\cS$ an equidimensional collection of compact subgroups thereof.
  \begin{enumerate}[(1)]
  \item The topologies $\cT_0$ and $\cT_1$ of \Cref{def:tpsbgp} coincide on $\cS$.
  \item If furthermore all members of $\cS$ have the same number of connected components then the topologies $\cT_{0,w}$, $\cT_0$ and $\cT_1$ all coincide on $\cS$.  \qedhere
  \end{enumerate}
\end{corollary}

In turn, that implies

\begin{corollary}\label{cor:alltopsame}
  Let $U$ be a Lie group and $G$ a compact group.

  On any collection of compact subgroups of $U$ isomorphic to $G$, the three topologies $\cT_{0,w}$, $\cT_0$ and $\cT_1$ coincide.  \qedhere
\end{corollary}

The discussion bifurcates naturally: we can consider
\begin{itemize}
\item spaces of (compact, hence also Lie) subgroups of $U$ which are all abstractly isomorphic;
\item and conditions under which sufficiently close compact subgroups of $U$ {\it must} be isomorphic.
\end{itemize}

The preceding discussion already provides criteria for the latter (automatic isomorphism).

\begin{proposition}\label{pr:1iso}
  If $G\le U$ is a compact subgroup of a Banach Lie group.
  \begin{enumerate}[(1)]
  \item If $K\le U$ is a compact subgroup sufficiently $\cT_{0,w}$-close to $G$ and of the same dimension, then the identity components of $K$ and $G$ are conjugate in $U$, and hence isomorphic.
  \item If $K$ is sufficiently $\cT_0$-close to $G$ and has the same dimension and the same number of components as $G$ then $K$ and $G$ are conjugate.
  \item\label{item:t1conj} The same conclusion holds if $K$ is sufficiently $\cT_1$-close to $G$.
  \end{enumerate}
\end{proposition}
\begin{proof}
  The argument is part of the proof of \Cref{pr:sametop}, namely of the implication \Cref{item:0s} $\Rightarrow$ \Cref{item:0w}.
\end{proof}

As for collections of groups we already know are mutually isomorphic, they have very rich structure:


\begin{theorem}\label{th:sameg}
  Let $U$ be a Lie group, $G$ a compact Lie group, and $\cS_{G\le U}$ the collection of compact subgroups of $U$ isomorphic to $G$, equipped with the topology $\cT=\cT_{\bullet}$ of \Cref{cor:alltopsame}.

  The map
  \begin{equation}\label{eq:xontos}
    \left(\text{continuous embeddings $G\le U$}\right)
    =:X\ni
    \varphi
    \xmapsto{\quad}
    \im(\varphi)
    \in \cS_{G\le U}
  \end{equation}
  is a principal $\Aut(G)$-fibration, where the automorphism group acts on the right by precomposition.

  In particular, the codomain $\cS_{G\le U}$ is an analytic manifold.
\end{theorem}
\begin{proof}
  Note first that $A:=\Aut(G)$ is a Lie group with compact identity component
  \begin{equation*}
    A_0\cong G_0/G_0\cap Z(G),
    \quad Z(G):=\text{center of }G
  \end{equation*}
  \cite[Corollary 6.69 (i) and sentence preceding Definition 6.70]{hm4}. 

  Next, the space $\Hom(G,U)$ of continuous morphisms $G\to U$ (equipped with the uniform topology) is an analytic manifold and the orbits of the action of $U$ by conjugation are open \cite[Corollary 2.4]{cp-us}, so
  \begin{equation*}
    X\subseteq \Hom(G,U)
  \end{equation*}
  is open and hence analytic. The fact that the precomposition action $X\times A\to X$ is also analytic is now an easy check.

  There are now a number of claims to verify. First, denote by $\cQ$ the quotient topology on $\cS:=\cS_{G\le U}$ induced by \Cref{eq:xontos}.
    
  \begin{enumerate}[(1)]
  \item {\bf \Cref{eq:xontos} is a principal $A$-fibration onto $(\cS,\cQ)$.} This is the familiar matter of equipping quotients with manifold structures. The conclusion will follow from \cite[\S III.1.5, Proposition 10]{bourb-lie-13}, once we have checked that its hypotheses hold:
    
    \begin{itemize}
    \item The action of $A$ on $X$ is {\it free}, i.e. non-trivial elements have no fixed points. Indeed, if $\varphi:G\to U$ is an embedding and $\alpha\in A$ a non-trivial automorphism of $G$ then $\varphi\ne \varphi\circ\alpha$.

    \item The action of $A$ on $X$ is {\it proper} \cite[\S III.4.1, D\'efinition 1]{bourb-top-1-4}, meaning that the map
      \begin{equation*}
        X\times A\ni(\varphi,\alpha)\xmapsto{\quad\cat{can}\quad}(\varphi,\varphi\circ\alpha)\in X\times X
      \end{equation*}
      is proper in the sense that $\cat{can}\times\id_Z$ is closed for any topological space $Z$ (i.e. maps closed sets to closed sets) \cite[\S I.10.1, D\'efinition 1]{bourb-top-1-4}.

      Given the freeness of the action, and hence the injectivity of $\cat{can}$, its properness follows from \cite[\S I.10.1, Proposition 2]{bourb-top-1-4} provided we show that it is closed.

      We have already noted that $A$ is a Lie group. On the other hand, $X$ is a uniformly topologized space of Lie-group valued maps on a compact space; because Lie groups are completely metrizable \cite[\S III.1.1, Proposition 1]{bourb-lie-13}, we can assess closure by means of sequence convergence. Consider, then, convergent sequences
      \begin{equation*}
        \varphi_n\xrightarrow[\quad n\quad]{}\varphi
        \quad\text{and}\quad
        \varphi_n\circ\alpha_n\xrightarrow[\quad n\quad]{} \varphi'
      \end{equation*}
      in $X$.

      Consider the restrictions of the automorphisms $\alpha_n\in A$ to the largest central torus $S:=Z_0(G_0)$ of the identity component $G_0$. If $S$ is $d$-dimensional then these can be identified with elements of
      \begin{equation*}
        GL_d:=GL_d(\bZ)\cong \Aut(\bT^d),\quad \bT^d = (\bS^1)^d = \text{the $d$-dimensional torus},
      \end{equation*}
      and the fact that $(\varphi_n)$ and $(\varphi_n\circ\alpha_n)$ are both uniformly convergent implies that said restriction $\alpha_n|_S$ take finitely many values in $GL_d$. By \Cref{le:cpctauto}, this means that $\alpha_n$ range over a {\it compact} subset of $A$, so by passing to a subsequence if necessary, we can assume that $\alpha_n\to \alpha\in A$. 

      We now know that convergent sequences in the image of $\cat{can}$ lift to convergent sequences in its domain, whence closure. 
      
    \item Finally, we have to check that for every embedding $\varphi:G\to U$ the orbit map
      \begin{equation}\label{eq:apa}
        A\ni\alpha\xmapsto{\quad}\varphi\circ\alpha\in X
      \end{equation}
      is an immersion at $1\in A$ \cite[\S III.1.5, Proposition 10 (b) and Proposition 9 (i)]{bourb-lie-13}.

      The tangent spaces $T_1A$ and $T_{\varphi}X$ both consist of continuous maps from $G$ to the tangent bundle of $U$ (see e.g. \cite[Theorem 2.2 (2)]{cp-us} for $T_{\varphi}X$), and the differential of \Cref{eq:apa} is nothing but composition with the differential $\varphi_{*}$. Because $\varphi$ is an embedding that differential is one-to-one, and the image of $T_1A$ in $T_{\varphi}X$ splits simply because it is finite-dimensional \cite[Theorem 7-3-6]{wil-tvs}.
      
    \end{itemize}

  \item {\bf The topologies $\cT$ and $\cQ$ coincide.} That $\cQ$ is at least as strong as $\cT$ is clear: uniformly close morphisms $G\to U$ have $\cT_0$-close images.

    Conversely, if $H$ is sufficiently $\cT_0$-close to $K$ in $\cS_{G\le U}$, the two are conjugate by some element $u\sim 1$ of $U$ by \Cref{th:closeconj}; but this also means that $H$ and $K$ are images of uniformly close morphisms
    \begin{equation*}
      G\cong H\le U
      \quad\text{and}\quad
      G\cong H\xrightarrow[\quad\cong\quad]{\Ad_u}K\le U,
    \end{equation*}
    hence $\cQ$-closeness. 
  \end{enumerate}  
\end{proof}

\Cref{th:sameg} has additional consequences on how well-behaved various classes of compact subgroups $G\le U$ are, even absent the mutual-isomorphism assumption.

\begin{theorem}\label{th:allsm}
  Let $U$ be a Lie group. Each of the following collections of compact subgroups of $U$ are analytic manifolds when equipped with the listed topologies:
  \begin{enumerate}[(1)]
  \item\label{item:allsm-d} the set $\cS_{d}$ of compact subgroups of any fixed dimension $d$, with the topology $\cT_{0}=\cT_1$;
  \item\label{item:allsm-dl} the set $\cS_{d,\ell}$ of compact subgroups of any fixed dimension $d$ with a fixed number $\ell$ of connected components, with $\cT_{0,w}=\cT_0=\cT_1$;
  \item\label{item:allsm-ss} the set $\cS_{ss}$ of semisimple (possibly disconnected) compact subgroups, with the topology $\cT_0=\cT_1$.
  \end{enumerate}
\end{theorem}
\begin{proof}
  That the topologies coincide as indicated in parts \Cref{item:allsm-d} and \Cref{item:allsm-dl} (on $\cS_{d}$ and $\cS_{d,\ell}$ respectively) follows from \Cref{cor:samedimsamecomp}.

  For any compact group $G\le U$ the compact subgroups of $U$ sufficiently $\cT_1$-close to $G$ are conjugate and in particular isomorphic to it (\Cref{pr:1iso}), so $\cS_d$ and $\cS_{d,\ell}$ are both partitioned by open subspaces $\cS_{G\le U}$ as in \Cref{th:sameg}, whence the analytic structure.

  As for item \Cref{item:allsm-ss}, observe that as a consequence of \Cref{pr:ssgood}, $\cT_0$-close compact {\it semisimple} subgroups of $U$  are isomorphic. It once more follows that we have a partition
  \begin{equation*}
    \cS_{ss} = \coprod_{G\text{ semisimple}}\cS_{G\le U}\text{ into open subsets};
  \end{equation*}
  conclude as before. 
\end{proof}


\section{Centralizer/normalizer continuity}\label{se:centnorm}

The title of the present section refers to the continuity of the operations of forming either centralizers or normalizers of compact subgroups $G\le U$ of a Banach Lie group, provided the relevant subgroup collections are topologized appropriately.

According to \Cref{pr:1iso} \Cref{item:t1conj}, sufficiently $\cT_1$-close compact subgroups of a Banach Lie groups are conjugate; it follows from this that so are their respective centralizers or normalizers. We are concerned here with the $\cT_0$ version of this remark. 

Centralizers and normalizers of compact groups are indeed Lie subgroups \cite[Corollary 2.22 and Proposition 2.26]{cp-us}, and the relevant topologies on the respective collections are variants of those of \Cref{def:tpsbgp}: because centralizers (or normalizers) will not, in general, be compact, actual convergence in any of the earlier topologies would be too much to expect; we rather need a notion of ``convergence in every bounded set'', which the following definition will be of use in formalizing.

\begin{definition}\label{def:adbdd}
  Let $U$ be a Banach Lie group with Lie algebra $\fu$.

  A subset $S\subseteq U$ is {\it ad-bounded} if
  \begin{equation*}
    \Ad(S):=\{\Ad_s\ |\ s\in S\}\subset \End(\fu)
  \end{equation*}
  is uniformly bounded for any Banach norm on $\fu$ compatible with its topological Lie algebra structure. The notion does not depend on the choice of norm.

  In general, for subsets $X$, $Y$ and $B$ of $U$ (of which $B$ will typically be ad-bounded), we write
  \begin{align*}
    X=_B Y&\text{ for }X\cap B = Y\cap B\text{ and}\\
    X\subseteq_B Y&\text{ for }X\cap B \subseteq Y\cap B.
  \end{align*}
  The same notation applies to other size-comparison symbols such as `$\supseteq$', `$\le$' (for inclusions of subgroups), etc.
\end{definition}

\begin{theorem}\label{th:cetnorm}
  Let $U$ be a Banach Lie group, $G\le U$ a compact subgroup, 
  \begin{equation*}
    G_n\xrightarrow[\quad n\quad]{0\text{ strong}}G
  \end{equation*}
  a 0-strongly-convergent sequence and $B\subseteq U$ an ad-bounded subset. 
  \begin{enumerate}[(1)]

  \item\label{item:centin} If $G_n\le G$, we have
    \begin{equation*}
      Z_U(G_n)  =_B Z_U(G)\text{ for large $n$.}
    \end{equation*}

  \item\label{item:centcj} In general,
    \begin{equation*}
      \Ad_u Z_U(G_n) =_B Z_U(G)
    \end{equation*}
    with $u\in U$ arbitrarily close to $1$ for large $n$.

  \item\label{item:normin} If $G_n\le G$, we have
    \begin{equation*}
      N_U(G_n)  \le_B N_U(G) \text{ for large $n$.}
    \end{equation*}

  \item\label{item:normcj} In general,
    \begin{equation*}
      \Ad_u N_U(G_n) \le_B N_U(G)
    \end{equation*}
    with $u\in U$ arbitrarily close to $1$ for large $n$.

  \end{enumerate}
\end{theorem}
\begin{proof}
  Parts \Cref{item:centcj} and \Cref{item:normcj} follow by \Cref{th:closeconj} from \Cref{item:centin} and \Cref{item:normin} respectively, so we focus on these latter claims. As far as \Cref{item:centin} is concerned, note that the inclusion
  \begin{equation*}
    Z_U(G_n)\ge Z_U(G)\Longrightarrow Z_U(G_n) \ge_B Z_U(G)
  \end{equation*}
  is in any case obvious: everything that centralizes $G$ also centralizes its subgroup $G_n$. For both \Cref{item:centin} and \Cref{item:normin}, then, the interesting inclusion is `$\le_B$'; in words, the claim is that for sufficiently large $n$, elements of $B$ that centralize (normalize) $G_n$ also centralize (respectively normalize) the larger group $G$.

  We write $\fg:=Lie(G)$ and $\langle-\rangle_{\bR}$ for the real span of a set. Note that
  \begin{itemize}
  \item for any small neighborhood $V_U$ of $1\in U$ large $n$ will ensure that 
    \begin{equation*}
      \left\langle \log\left(G_n\cap V_U\right)\right\rangle_{\bR}
      =
      \left\langle \log\left(G\cap V_U\right)\right\rangle_{\bR}
      =
      \fg;
    \end{equation*}
  \item and hence 
    \begin{equation}\label{eq:ggnexp}
      G = G_n\cdot \exp \left\langle \log\left(G_n\cap V_U\right)\right\rangle_{\bR} = G_n\cdot \exp\fg = G_n\cdot G_0
    \end{equation}
    for large $n$, given that the exponential of a compact connected Lie group (such as the identity component $G_0$) is onto \cite[Theorem 9.49 (ii)]{hm4}.
  \end{itemize}
  We can select nested neighborhoods $V'_U\subset V_U$ so that 
  \begin{equation}\label{eq:adubb}
    \Ad_u V'_{U} \subset V_{U},\ \forall u\in B
  \end{equation}
  (this is where the ad-boundedness of $B$ is used).
  

  
  An operator $\Ad_u$ that centralizes (normalizes) $G_n$ will then also centralize $\log(G_n\cap V'_U)$ or, respectively, map it to a subset of $\log(G_n\cap V_U)$ still spanning $\fg$, whence the conclusion by \Cref{eq:ggnexp}.
\end{proof}

Even in the finite-dimensional case, one cannot hope to do better, in \Cref{th:cetnorm}, than equality/containment ``along bounded sets'':

\begin{example}\label{ex:onlyalongb}
  Let $U$ be the semidirect product $\bT^2\rtimes \bZ$, where the generator of $\bZ$ acts on
  \begin{equation*}
    \bT^2\cong \bR^2/\bZ^2
    \text{ via }
    \begin{pmatrix}
      1&1\\
      0&1
    \end{pmatrix}.
  \end{equation*}
  The subgroup $G\le U$ will be the diagonal torus
  \begin{equation}\label{eq:gdiag}
    G := \{(z,z)\ |\ z\in \bS^1\}\subset \bT^2,
  \end{equation}
  while $G_n\subset G$ is the finite subgroup of $2^n$th roots of unity therein. The 0-strong convergence $G_n\to G$ is clear, but while for every $n$ the subgroup
  \begin{equation*}
    2^n \bZ\subset \bZ\subset U=\bT^2\rtimes \bZ
  \end{equation*}
  centralizes $G_n$, none of these groups normalize $G$. For that reason, neither \Cref{item:centin} nor \Cref{item:normin} holds without the `$B$' subscript on the equality/inclusion symbols.
\end{example}

\begin{remark}\label{re:notlin}
  Note in passing that the group $U\cong \bT^2\rtimes \bZ$ of \Cref{ex:onlyalongb} is not linear in the sense of \cite[Definition 5.32]{hm4}:
  \begin{itemize}
  \item On the one hand, the proof of \cite[Proposition 5.54 (ii)]{hm4} argues (in the portion of the proof appealing to \cite[Lemma 5.15 (i)]{hm4}) that if $g\in U$ is sufficiently close to $1$, an element which centralizes it also centralizes the 1-parameter group containing $\log g$.
  \item On the other, this is not the case in \Cref{ex:onlyalongb}, since there are $2^n$th roots of unity arbitrarily close to $1$, each centralized by some $2^n\bZ$ that fails to centralize the diagonal \Cref{eq:gdiag}.
  \end{itemize}
\end{remark}

\section{A supplement: automatic compactness}\label{se:autocpct}

The constant theme throughout the discussion below is that of Lie subgroups $K\le U$ of (a Banach Lie group) $U$ contained in small neighborhoods of a fixed compact group $G$. As indicated above, we will be interested in {\it compact} $K$, so it is perhaps worth pausing to note that compactness is automatic under the appropriate assumptions.

\begin{theorem}\label{th:autocpct}
  Let $U$ be a Banach Lie group and $G\le U$ a compact subgroup.

  There is a neighborhood $V\supseteq G$ such that every Lie subgroup $K\le U$ contained in $V$ is compact.
\end{theorem}

It follows that throughout the preceding discussion {\it compact} (hence Lie) subgroups of $U$ sufficiently $\cT_{0,w}$-close to a compact $G\le U$ may as well, to begin with, have been assumed only Lie.

The proof requires some preliminary work. 
And an additional piece of notation: having fixed a complete metric $d$ topologizing $U$,
\begin{equation*}
  G_{\varepsilon}:=\{x\in U\ |\ d(x,G)<\varepsilon\}
\end{equation*}
denotes the $\varepsilon$-neighborhood of $G$. We can work with these in place of the generic $V$ of \Cref{th:autocpct}.

\begin{lemma}\label{le:findim}
  Let $G\le U$ be a compact subgroup of a Banach Lie group. 
  
  For small $\varepsilon$, a Lie subgroup $K\subset G_{\varepsilon}$ is at most $(\dim G)$-dimensional.
\end{lemma}
\begin{proof}
  Consider a sufficiently small origin neighborhood $B\subset \fu:=Lie(U)$ so that in particular $\exp:\fu\to U$ induces mutually-inverse isomorphisms
  \begin{equation*}
    \begin{tikzpicture}[auto,baseline=(current  bounding  box.center)]
      \path[anchor=base] 
      (0,0) node (l) {$B$}
      +(4,0) node (r) {$\exp(B)$.}
      ;
      \draw[->] (l) to[bend left=10] node[pos=.5,auto] {$\scriptstyle \exp$} (r);
      \draw[->] (r) to[bend left=10] node[pos=.5,auto] {$\scriptstyle \log$} (l);
    \end{tikzpicture}
  \end{equation*}
  The almost-containment of $K\cap \exp(B)$ in $G\cap \exp(B)$ implies that
  \begin{equation*}
    \fk\cap B\text{ is almost contained in }\fg\cap B,\text{ where }\fk:=Lie(K),\ \fg:=Lie(G). 
  \end{equation*}
  This implies in particular that if
  \begin{equation*}
    K\le G_{\varepsilon},\ \text{small }\varepsilon>0
  \end{equation*}
  then $K$ is finite-dimensional, with $\dim K\le \dim G$.
\end{proof}

\begin{lemma}\label{le:nor}
  Let $G\le U$ be a compact subgroup of a Banach Lie group.

  For small $\varepsilon$, a Lie subgroup $K\subset G_{\varepsilon}$ has compact identity component $K_0$.
\end{lemma}
\begin{proof}
  If not, it contains a copy of $(\bR,+)$ as a closed subgroup (by basic Lie-group structure theory; e.g. \cite[\S 4.13, first stated Theorem]{mz}). Being a 1-parameter subgroup of $K$, that copy must be of the form
  \begin{equation*}
    \{\exp(tv)\ |\ t\in \bR\}\text{ for some }v\in \fk. 
  \end{equation*}
  Because all of $\exp\left(\bR v\right)$ is contained in $G_{\varepsilon}$, there will be large $t_0$ with $w:=\log\exp(t_0v)\in B$. The 1-parameter group $\exp\left(\bR w\right)$ then commutes with $\exp\left(\bR w\right)$, and the two generate a 2-dimensional integral subgroup \cite[\S III.6.2, Definition 1 and Theorem 2 (i)]{bourb-lie-13} $L\le U$, with Lie algebra $\fl\cong \bR^2$.

  Being abelian and non-compact, $L$ is (isomorphic to) one of $\bR^2$ or $\bR\times \bS^1$ \cite[Theorem 4.2.4]{de}. Neither of these, though, contains
  \begin{itemize}
  \item closed 1-parameter subgroups (such as $\exp\left(\bR v\right)$);
  \item with non-trivial elements (such as $\exp(t_0v)$);
  \item that also belong to a {\it second} 1-parameter subgroup (such as $\exp\left(\bR v\right)$).
  \end{itemize}
  The contradiction shows that $K_0$ is compact.
\end{proof}

For the sequel, note the alternative possible phrasing of \Cref{th:autocpct} (as in \Cref{re:altform} \Cref{item:seq}): if
\begin{equation}\label{eq:kntog}
  K_n\xrightarrow[n]{0\text{ weak}}G
\end{equation}
is sequence of Lie subgroups of $U$ 0-weakly converging to a compact subgroup $G$, then the $K_n$ are eventually compact. We cast the problem in this manner whenever convenient.

We handle sequences of {\it discrete} subgroups separately.

\begin{proposition}\label{pr:discfinite}
  Let $G\le U$ be a compact subgroup of a Banach Lie group and \Cref{eq:kntog} a sequence of discrete closed subgroups of $U$.

  For small $\varepsilon$, every discrete closed subgroup $K\subset G_{\varepsilon}$ is finite.
\end{proposition}
\begin{proof}
  There are a few stages to the argument.

  \begin{enumerate}[(I)]
  \item\label{item:mustbetorsion} {\bf For small $\varepsilon$, $K\subseteq G_{\varepsilon}$ is torsion.}

    Suppose, on the contrary, that we can find infinite-order elements $x\in U$ generating closed copies
    \begin{equation*}
      \{x^k\ |\ k\in \bZ\}\subseteq G_{\varepsilon}
    \end{equation*}
    of $(\bZ,+)$ for arbitrarily small $\varepsilon$. Each $x^k$ is $\varepsilon$-close to some $g_{k}\in G$; because $G$ is compact, some subsequence of $(g_k)_k$ converges to some $g\in G$, so there will be pairs
    \begin{equation*}
      d(x^{k_1},x^{k_2})<2\varepsilon
      \quad\text{with}\quad
      k_1<k_2
      \quad\text{and}\quad
      k_2-k_1
      \text{ arbitrarily large}.
    \end{equation*}
    It follows, then, that we can find infinite-order $x$ arbitrarily close to $1$. But then the 1-parameter group
    \begin{equation*}
      \left\{\exp\left(t\log x\right)\ |\ t\in \bR\right\}\cong (\bR,+)
    \end{equation*}
    is closed and, being a union of $K$-translates of the segment
    \begin{equation*}
      \exp\left([0,1]\log x\right),
    \end{equation*}
    can be chosen arbitrarily $\cT_{0,w}$-close to $G$. The contradiction now follows from \Cref{le:nor}.

  \item\label{item:torsion} {\bf Torsion groups.} More explicitly, the goal is to show that for small $\varepsilon$, torsion discrete $K\le G_{\varepsilon}$ is finite.

    For any $\varepsilon>0$, we can arrange for a sufficiently small $\varepsilon'>0$ so that all $K\subseteq G_{\varepsilon'}$ are covered by at most $N$ $K$-translates of the $\varepsilon$-neighborhood
    \begin{equation*}
      U_{1,\varepsilon}:=\{u\in U\ |\ d(1,u)<\varepsilon\}\subset U. 
    \end{equation*}

    For sufficiently small $\varepsilon$ the iterated commutators $U_{1,\varepsilon}^{(n)}$ decrease to $\{1\}$ by \Cref{le:iteratedcomm}, so every group generated by finitely many elements of $K_{1,\varepsilon}:= K\cap B_{1,\varepsilon}$ is nilpotent.

    Consider a group $L\le K$ as in the preceding paragraph, generated by finitely many $x\in K_{1,\varepsilon}$. We are assuming $x$ to be torsion; if $\varepsilon>0$ is small enough for `$\log$' to make sense, the 1-parameter groups
    \begin{equation*}
      {}_x\bS^1:=\left\{\exp(t\log x)\ |\ t\in \bR\right\}
    \end{equation*}
    are circles, and they generate a compact, connected, nilpotent subgroup of $U$: nilpotence follows inductively on the nilpotent length of $L$ from the fact that elements centralizing $x$ also centralize the line $\bR \log x\le Lie(U)$, etc. Because compact, connected, nilpotent Lie groups are abelian and hence tori (consequence of \cite[Theorem 6.19]{hm4}, say), $L$ is in fact abelian. In conclusion, all elements of $K$ sufficiently close to $1$ commute.

    Note furthermore that the tori in which finitely-generated subgroups $L\le K$ in the preceding discussion are embedded will be arbitrarily $\cT_{0,w}$-close to $G$ by the same argument as in the proof of \Cref{le:nor}, so by \Cref{le:findim} they cannot increase in dimension arbitrarily. It follows that the torsion abelian group generated by $K_{1,\varepsilon}$ is finite (because there is a cap on the dimension of a torus containing it), so $K_{1,\varepsilon}$ itself is finite.

    $K$ being a union of $\le N$ translates of its finite subset $K_{1,\varepsilon}$, it must be finite.
  \end{enumerate}
  Since \Cref{item:mustbetorsion} reduces the problem to torsion groups and \Cref{item:torsion} disposes of those, we are done.
\end{proof}

Some notation pertinent to \Cref{le:iteratedcomm}: for subsets $X,Y\subseteq U$ write
\begin{equation*}
  [X,Y] := \{[x,y]:=xyx^{-1}y^{-1}\ |\ x\in X,\ y\in Y\}
\end{equation*}
and
\begin{equation*}
  X^{(0)}:=X,\ X^{(n+1)}:=[X,X^{(n)}],\ \forall n\ge 0
\end{equation*}
for its iterated commutators (e.g. \cite[Proposition 8.8]{ragh}). 

The following result is worked out in the finite-dimensional case in \cite[Theorem 8.8 and Corollary 8.15]{ragh}; the general (infinite-dimensional) version follows essentially as sketched in \cite[proof of Proposition 4.9.2]{thurst-3mfld-unpublished}: the commutator map
\begin{equation*}
  U\times U\xrightarrow[]{\quad [-,-]\quad} U
\end{equation*}
has vanishing derivative at $(1,1)$, so forming commutators with sufficiently small neighborhoods $V\ni 1$ in $U$ scales distances in the Lie algebra $Lie(U)$ by at most $C$ for some $0<C<1$, etc.

\begin{lemma}\label{le:iteratedcomm}
  If $V$ is a sufficiently small neighborhood of $1\in U$ in a Banach Lie group, the sequence of iterated commutators $V^{(n)}$ is decreasing by inclusion and 0-converges to $\{1\}$. \qedhere
\end{lemma}

\pf{th:autocpct}
\begin{th:autocpct}
  Let \Cref{eq:kntog} be a 0-weakly-convergent sequence of Lie subgroups. By \Cref{pr:discfinite}, it will be enough to reduce the problem to the case of {\it discrete} $K_n$.
  
  \Cref{th:closeconj} and the compactness of the identity components $K_{n,0}$ (\Cref{le:nor}) ensure that
  \begin{equation*}
    u_n\cdot K_{n,0}\cdot u_{n}^{-1}\le G \text{ for some } u_n\to 1\in U. 
  \end{equation*}
  passing to a subsequence if needed, so that the conjugates $u_n\cdot K_{n,0}\cdot u_{n}^{-1}$ 0-(strongly-)converge to $H\le G$, we can assume that all $K_{n,0}$ are {\it contained} in $H$, and furthermore, by \Cref{le:w00}, that
  \begin{equation*}
    K_{n,0}\xrightarrow[\quad n\quad]{0\text{ strong}} K
  \end{equation*}
  for some closed (hence compact) connected $K\le G$.

  Now, is contained in the $\varepsilon$-neighborhood
  \begin{equation*}
    G_{\varepsilon} = G\cdot U_{1,\varepsilon},\quad U_{1,\varepsilon}:=\{u\in U\ |\ d(1,u)<\varepsilon\}
  \end{equation*}
  of the compact group $G$ for large $n$, meaning that it is ad-bounded in the sense of \Cref{def:adbdd}: compact subsets of $U$ are ad-bounded, as are sufficiently small neighborhoods of $1\in U$, and a product
  \begin{equation*}
    X\cdot Y:=\{xy\ |\ x\in X,\ y\in Y\}\subseteq U
  \end{equation*}
  of two ad-bounded sets is again ad-bounded.
 
  Because $K_n$ normalizes $K_{n,0}$, for large $n$ it will also normalize $K$ by \Cref{th:cetnorm} \Cref{item:normin}. Relabeling $K\cdot K_n$ as $K_n$, we now have $K_{n,0}=K$ for all $n$ and hence all $K_n$ are contained in the normalizer $N_U(K)$ (a Lie subgroup of $U$ by \cite[Corollary 2.6]{cp-us}).

  Now substitute
  \begin{itemize}
  \item $N_U(K)/K$ for $U$;
  \item $N_G(K)/K$ for $G$;
  \item and the discrete subgroups
    \begin{equation*}
      K_n/K\le N_U(K)/K\text{ for }K_n\text{ respectively}.
    \end{equation*}
  \end{itemize}
  As explained, an application of \Cref{pr:discfinite} finishes the proof of \Cref{th:autocpct}.
\end{th:autocpct}

\addcontentsline{toc}{section}{References}

\Addresses

\end{document}